\documentclass{nlaauth}
\usepackage{mathptmx}
\usepackage{amsfonts,amssymb,amsbsy,amsmath}
\usepackage{pifont,latexsym,ifthen,theorem,rotating,calc}
\allowdisplaybreaks
\numberwithin{equation}{section}

\newtheorem{theorem}{Theorem}[section]

\newtheorem{example}{Example}[section]

\newtheorem{remark}{Remark}[section]

\begin{document}
\NLA{1}{23}{1}{00}{14}

\runningheads{H.Y. LI AND Y.M. WEI} {Rigorous perturbation bounds for LU and QR factorizations}

\title{New rigorous perturbation bounds for the LU and QR factorizations}

\author{Hanyu Li\affilnum{1}\corrauth and Yimin Wei\affilnum{2}}

\address{\affilnum{1}College of Mathematics and Statistics, Chongqing University,
Chongqing, 401331, P. R. China.\\\affilnum{2}School of Mathematical Sciences and Key Laboratory
of Mathematics for Nonlinear Sciences,
Fudan University, Shanghai,
200433, P. R. China.}

\corraddr{Hanyu Li, College of Mathematics and Statistics, Chongqing University,
Chongqing, 401331, P. R. China. E-mail address: lihy.hy@gmail.com or hyli@cqu.edu.cn (H.Y. Li); yimin.wei@gmail.com or ymwei@fudan.edu.cn (Y.M. Wei).}


\cgsn{The work is supported by the National Natural Science Foundation of China}{11201507, 11271084}

\received{} \revised{} \noaccepted{}

\begin{abstract}
Combining the modified matrix-vector equation approach with the
technique of Lyapunov majorant function and the Banach fixed point
principle, we obtain new rigorous perturbation bounds for the LU and
QR factorizations with normwise or componentwise perturbations in
the given matrix, where the componentwise perturbations have the
form of backward errors resulting from the standard factorization
algorithms. Each of the new rigorous perturbation bounds is a
rigorous version of the first-order perturbation bound derived by
the matrix-vector equation approach in the literature, and we
present their explicit expressions. These bounds improve the results
given by Chang and Stehl\'{e} [{\em SIAM Journal on Matrix Analysis and Applications} 2010; 31:2841--2859]. Moreover, we derive  new tighter
first-order perturbation bounds including two optimal ones for the
LU factorization, and provide the explicit expressions of the
optimal first-order perturbation bounds for the LU and QR
factorizations.
\end{abstract}



\section{INTRODUCTION}
Let $\mathbb{R}^{m \times n} $ be the set of $m \times n$ real
matrices and $\mathbb{R}^{m\times n}_r $ be the subset of
$\mathbb{R}^{m \times n} $ with rank $r$. Let $I_r $ be the identity
matrix of order $r$ and $A^T$ be the transpose of the matrix $A$.

For a matrix $A\in \mathbb{R}^{n\times n}$, if its leading principal
sub-matrices are all nonsingular, then there exists a unique unit
lower triangular matrix $L\in \mathbb{R}^{n\times n}$ and a unique
upper triangular matrix $U\in \mathbb{R}^{n\times n}$ such that
\begin{eqnarray}
A=LU.
\end{eqnarray}
The factorization is called the LU factorization of the matrix $A$,
and the matrices $L$ and $U$ are referred to as the LU factors. The
LU factorization is a basic and effective tool in numerical linear
algebra (see e.g., \cite{Golub12,Higham02}).

For a matrix $A\in \mathbb{R}^{m\times n}_n$, there exists a unique
matrix $Q\in \mathbb{R}^{m\times n}$ with orthonormal columns, i.e.,
$Q^TQ=I_n$, and a unique upper triangular matrix $R\in
\mathbb{R}^{n\times n}$ with positive diagonal elements such that
\begin{eqnarray}
A=QR.
\end{eqnarray}
The factorization is called the QR factorization of the matrix $A$,
and the matrices $Q$ and $R$ are named after the orthonormal factor
and the triangular factor, respectively. The QR factorization is an
important tool in matrix computations (see e.g., \cite{Golub12,Higham02}).

For the LU and QR factorizations, their applications, algorithms,
and stability of algorithms have been considered (see e.g., \cite{Golub12,Higham02,Bai99}). Since the object matrix $A$ may be
contaminated by the errors from measurement, modeling, and so on,
and the numerical algorithms will introduce rounding errors in
computing these factorizations, the computed factors may not be the
exact ones. Naturally, it is important to know how much the factors
may change when the original matrix changes. Therefore, several
scholars discussed the perturbation analysis of the LU and QR
factorizations. The first rigorous perturbation bounds for the LU
factorization was derived by Barrlund \cite{Barr91} when the original matrix
has the normwise perturbation. Here, a bound is said to be {\it
rigorous} if it doesn't neglect any higher-order terms. Later, using
a different approach, Stewart \cite{Stewart93} presented the first-order
perturbation bounds. These results were improved in
\cite{Stewart97}. For the QR factorization, the first rigorous
perturbation bounds with normwise perturbation were given in \cite{Stewart77},
which were further modified and improved by
Sun \cite{Sun91}.  Sun \cite{Sun91} also provided the first-order
perturbation bounds, which were obtained by Stewart \cite{Stewart93} too
using a different approach. Later, Sun \cite{Sun95} presented new
rigorous perturbation bounds for the orthonormal factor $Q$ alone,
from which an improved first-order perturbation bound was derived.
This bound was also given in \cite{Bha94}.

In 1996, Chang {\it et al}. \cite{Chang96} proposed the refined matrix equation
approach and the matrix-vector equation approach, which can be used
to apply the first-order perturbation analysis of many matrix
factorizations, such as, the Cholesky, LU, QR, and SR factorizations \cite{Chang96,Chang97a,Chang97b,Chang98,Chang02,Chang98a,Chang98c,Chang01,Chang97},
when the original matrix has normwise or componentwise
perturbations. Here, the componentwise perturbation have the form
of backward errors for the standard factorization algorithms (see e.g.,
\cite{Bai99}). This class of perturbations was first investigated by
Zha \cite{Zha93} for the QR factorization. The new first-order
perturbation bounds with the above two approaches improve the previous
ones greatly. Recently, a new approach, the combination of the
classic and refined matrix equation approaches, was provided by
Chang {\it et al}. to study the rigorous perturbation bounds for
some matrix factorizations \cite{Chang12,Chang11,Chang10,Chang12a}. With their
approach, the new rigorous perturbation bounds can be much smaller
than the previous ones derived by the classic matrix equation
approach. In addition, the rigorous perturbation bounds for the
Cholesky factorization can also be obtained by combining the
matrix-vector equation approach and the results in \cite[Theorem
3.1]{Stewart73}; the reader can refer to \cite{Chang96} or
\cite{Chang97a}. These bounds are tighter than the ones in
\cite{Chang10}. However, the above technique can not be applied to
the LU factorization. The main reason is that Theorem 3.1 in
\cite{Stewart73} can not be used any longer. Furthermore, the
rigorous bounds derived by the above technique have no explicit
expressions and then it is difficult to interpret and understand
them.

In this paper, we combine the modified matrix-vector equation
approach, the technique of Lyapunov majorant function (see,
e.g., \cite[Chapter 5]{Konstantinov03}), and the Banach fixed point
principle (see,
e.g., \cite[Appendix D]{Konstantinov03}) to
investigate the rigorous perturbation bounds for the LU
factorization. Moreover, the rigorous perturbation bounds for the
triangular factor $R$ of the QR factorization are also obtained by
using the above approach. The new bounds for the LU and  QR
factorizations can be regarded as the rigorous versions of the
first-order perturbation bounds derived by the matrix-vector
equation approach in \cite{Chang97a}, \cite{Chang98a},
\cite{Chang01}, and \cite{Chang97}, have the explicit expressions,
and improve the corresponding ones in \cite{Chang10} and
\cite{Chang12a}.

The rest of this paper is organized as follows. Section 2 presents
some notation and preliminaries. The rigorous perturbation bounds
for the LU and QR factorizations with normwise or componentwise
perturbations are given in Sections 3 and 4, respectively. In
particular, new tighter first-order perturbation bounds for the LU
factorization and the explicit expressions of the optimal first-order perturbation bounds for the LU and QR
factorizations are also provided in these two sections. Finally, we present the
concluding remarks of the whole paper.

\section{NOTATION AND PRELIMINARIES}
Given the matrix $A=(a_{ij}) \in \mathbb{R}^{m \times n}$, the
symbols $A^{\dag}$, $\left\| A \right\|_2 $, and $\left\| A
\right\|_F $ stand for its Moore-Penrose inverse (see,
e.g., \cite[Chapter III]{Stewart90}), spectral
norm, and Frobenius norm, respectively,
$\kappa_2(A)=\left\|A^{\dag}\right\|_2\left\|A\right\|_2$ denotes
its condition number, and $|A|$ is defined by $|A|=(|a_{ij}|)$. For
the above two norms, the following relations hold (see,
e.g., \cite[page
80]{Stewart90}),
\begin{eqnarray}
{\left\| {XYZ} \right\|_F} \leq {\left\| X \right\|_2}{\left\| Y \right\|_F}{\left\| Z \right\|_2},\quad {\left\| {XYZ} \right\|_2} \leq
{\left\| X \right\|_2}{\left\| Y \right\|_2}{\left\| Z \right\|_2},
\end{eqnarray}
whenever the matrix product $XYZ$ is defined. Note that the
Frobenius norm is monotone (see,
e.g., \cite[Chapter 6]{Higham02}). That is, for a matrix
$B=(b_{ij}) \in \mathbb{R}^{m \times n}$, if $|A|\leq |B|$, then
$\left\| A \right\|_F=\left\| |A| \right\|_F\leq \left\| |B
|\right\|_F =\left\| B \right\|_F$. Here $A\leq B$ means $a_{ij}\leq
b_{ij}$ for each $i=1,2,\ldots,m$, $j=1,2,\ldots,n$. In addition,
for a matrix 2-tuple $C=\left[ {\begin{array}{*{20}c}
   A   \\
   B  \\
\end{array}} \right]$, we define the `generalized matrix norm' (see,
e.g., \cite[page
13]{Konstantinov03}) by
\begin{eqnarray}
|||C| ||=\left[ {\begin{array}{*{20}c}
   \left\| A \right\|_F   \\
   \left\| B \right\|_F  \\
\end{array}} \right].
\end{eqnarray}

For the matrix $A=[a_1,a_2,\cdots, a_n]=(a_{ij})\in
\mathbb{R}^{n\times n}$, we denote the vector of the first $i$
elements of $a_j$ by $a^{(i)}_j$ and the vector of the last $i$
elements of $a_j$ by $a^{[i]}_j$. With these, we adopt the operators
as in \cite{Chang97a},
\begin{eqnarray*}
&{\rm {uvec}}(A): = \left[ {\begin{array}{*{20}{c}}
   {a_1^{(1)}}  \\
   {a_2^{(2)}}  \\
    \vdots   \\
   {a_n^{(n)}}  \\
\end{array}} \right]\in\mathbb{R}^{\nu_1},\ {\rm {slvec}}(A): = \left[ {\begin{array}{*{20}{c}}
   {a_1^{[n-1]}}  \\
   {a_2^{[n-2]}}  \\
    \vdots   \\
   {a_{n-1}^{[1]}}  \\
\end{array}} \right]\in
\mathbb{R}^{\nu_2},\ {\rm {vec}}(A): = \left[ {\begin{array}{*{20}{c}}
   {a_1}  \\
   {a_2}  \\
    \vdots   \\
   {a_n}  \\
\end{array}} \right]\in
\mathbb{R}^{n^2},\\
&{\rm {up}}(A): = \left[ {\begin{array}{*{20}{c}}
   {\frac{1}{2}a_{11}} & {{a_{12}}} &  \cdots  & {{a_{1n}}}  \\
   0 & {\frac{1}{2}a_{22}}&  \cdots  & {{a_{2n}}}  \\
    \vdots  &  \vdots  &  \ddots  &  \vdots   \\
   0 & 0 &  \cdots  & {\frac{1}{2}a_{nn}}  \\
\end{array}} \right]\in
\mathbb{U}_{n},\ {\rm {ut}}(A): = \left[ {\begin{array}{*{20}{c}}
   {a_{11}} & {{a_{12}}} &  \cdots  & {{a_{1n}}}  \\
   0 & {a_{22}}&  \cdots  & {{a_{2n}}}  \\
    \vdots  &  \vdots  &  \ddots  &  \vdots   \\
   0 & 0 &  \cdots  & {a_{nn}}  \\
\end{array}} \right]\in
\mathbb{U}_{n},
\end{eqnarray*}
and
\begin{eqnarray*}
& {\rm {slt}}(A): =A-{\rm {ut}}(A)\in
\mathbb{SL}_{n},
\end{eqnarray*}
where $\nu_1=n(n+1)/2$, $\nu_2=n(n-1)/2$, and $\mathbb{U}_{n}$ and
$\mathbb{SL}_{n}$ denote the sets of $n \times n$ real upper
triangular and strictly lower triangular matrices, respectively. Considering the structures of
these operators, we have
\begin{eqnarray}
{\rm {uvec}}(A)=M_{{\rm {uvec}}}{\rm {vec}}(A),\ {\rm {slvec}}(A) =M_{{\rm {slvec}}}{\rm {vec}}(A),
\end{eqnarray}
and
\begin{eqnarray}
{\rm {vec}}({\rm {up}}(A))=M_{{\rm {up}}}{\rm
{vec}}(A),\ {\rm {vec}}({\rm {ut}}(A))=M_{{\rm {ut}}}{\rm
{vec}}(A),\ {\rm {vec}}({\rm {slt}}(A))=M_{{\rm {slt}}}{\rm
{vec}}(A),
\end{eqnarray}
where
\begin{eqnarray*}
&&M_{{\rm {uvec}}}={\rm diag}\left(J_1, J_2,\cdots,J_n\right)\in \mathbb{R}^{\nu_1\times n^2},\ J_i=\left[I_i,\ 0_{i\times (n-i)}\right]\in \mathbb{R}^{i\times n},\\
&&M_{{\rm {slvec}}}=\left[{\rm diag}\left(\widetilde{J}_1,\widetilde{J}_2,\cdots,\widetilde{J}_{n-1}\right),\ 0_{\nu_2\times n}\right]\in \mathbb{R}^{\nu_2\times n^2},\ \widetilde{J}_i
=\left[0_{(n-i)\times i},\ I_{n-i}\right]\in \mathbb{R}^{(n-i)\times n},\\
&&M_{{\rm {up}}}={\rm diag}\left(S_1,S_2,\cdots,S_n\right)\in \mathbb{R}^{n^2\times n^2},\ S_i={\rm diag}\left(I_{i-1},1/2,0_{(n-i)\times (n-i)}\right)\in \mathbb{R}^{n\times n},\\
&&M_{{\rm {ut}}}={\rm diag}\left(\widetilde{S}_1,\widetilde{S}_2,\cdots,\widetilde{S}_n\right)\in \mathbb{R}^{n^2\times n^2},\ \widetilde{S}_i={\rm diag}\left(I_i,0_{(n-i)\times (n-i)}\right)\in \mathbb{R}^{n\times n},\\
&&M_{{\rm {slt}}}={\rm
diag}\left(\widehat{S}_1,\widehat{S}_2,\cdots,\widehat{S}_{n-1},0_{n\times
n}\right)\in \mathbb{R}^{n^2\times n^2},\ \widehat{S}_i={\rm
diag}(0_{i\times i},I_{n-i})\in \mathbb{R}^{n\times n}.
\end{eqnarray*}
Here, $0_{s\times t}$ is the $s\times t$ zero matrix. It is easy to verify that
\begin{eqnarray}
M_{{\rm {uvec}}}M_{{\rm {uvec}}}^T=I_{\nu_1},\ M_{{\rm {slvec}}}M_{{\rm {slvec}}}^T=I_{\nu_2},
\end{eqnarray}
and
\begin{eqnarray}
M_{{\rm {uvec}}}^TM_{{\rm {uvec}}}=M_{{\rm {ut}}},\ M_{{\rm
{slvec}}}^TM_{{\rm {slvec}}}=M_{{\rm {slt}}}.
\end{eqnarray}
Let ${\rm {uvec^\dag}}:\mathbb{R}^{\nu_1}\rightarrow
\mathbb{R}^{n\times n}$ be the right inverse of the operator `uvec'
such that ${\rm {uvec}}\cdot{\rm {uvec^\dag}}=1_{\nu_1\times \nu_1}$
and ${\rm {uvec^\dag}}\cdot{\rm {uvec}}={\rm {ut}}$. Then the matrix
of the operator `${\rm {uvec^\dag}}$' is $M_{{\rm {uvec}}}^T$.
That is, ${\rm {uvec^\dag}}(A)=M^T_{{\rm {uvec}}}{\rm {vec}}(A)$. Similarly, we can
define the right inverse of the operator `slvec' by `${\rm
{slvec^\dag}}$', whose matrix is $M_{{\rm {slvec}}}^T$. Some results mentioned
above can be found in \cite{Konstantinov02}.

Let $A =(a_{ij})\in {\mathbb{R}^{m \times n}}$ and $B \in
{\mathbb{R}^{p \times q}}$. The {\em Kronecker product} between $A$ and $B$ is defined
by (see,
e.g., \cite[Chapter 4]{Horn91}),
\[A \otimes B = \left[ {\begin{array}{*{20}{c}}
{{a_{11}}B}&{{a_{12}}B}& \cdots &{{a_{1n}}B}\\
{{a_{21}}B}&{{a_{22}}B}& \cdots &{{a_{2n}}B}\\
 \vdots & \vdots & \ddots & \vdots \\
{{a_{m1}}B}&{{a_{m2}}B}& \cdots &{{a_{mn}}B}
\end{array}} \right]\in {\mathbb{R}^{mp \times nq}}.\]
It follows from \cite[Chapter 4]{Horn91} that
\begin{eqnarray}
{\rm vec}(AXB) = \left({B^T} \otimes A\right){\rm vec}(X)
\end{eqnarray}
and
\begin{eqnarray}
\Pi_{mn} {\rm vec}(A) = {\rm vec}({A^T}),
\end{eqnarray}
where $X\in {\mathbb{R}^{n \times p}}$, and $\Pi_{mn} \in {\mathbb{R}^{mn \times mn}}$ is called the vec-permutation matrix and
can be expressed explicitly by
$$\Pi_{mn}  = \sum\limits_{i = 1}^n {\sum\limits_{j = 1}^m {{E_{ij}}(m \times n)} }  \otimes {E_{ji}}(n \times m).$$
In the above expression, ${E_{ij}}(m \times n) =
e_i{(m)}{(e_j{(n)})^T} \in {\mathbb{R}^{m \times n}}$ denotes the
$(i,j)$-th elementary matrix and $e_i{(m)}$ is the vector ${\left[
{0,0, \cdots ,0,1,0, 0,\cdots ,0} \right]^T} \in {\mathbb{R}^m}$,
i.e., the $1$ in the $i$-th component. In addition, from \cite[Chapter 4]{Horn91}, we also have that if $A$ and $B$ are nonsingular, then $A \otimes B$ is also nonsingular and
\begin{eqnarray}
(A \otimes B)^{-1}=A^{-1} \otimes B^{-1}.
\end{eqnarray}

\section{PERTURBATION BOUNDS FOR THE LU FACTORIZATION}
Assume that the matrices $A$, $L$, and $U$ in (1.1) are perturbed as
\begin{eqnarray*}
A\rightarrow A+\Delta A,\ L\rightarrow L+\Delta L,\ U\rightarrow U+\Delta U,
\end{eqnarray*}
where $\Delta A \in {\mathbb{R}^{n \times n}}$, $\Delta L\in \mathbb{SL}_{n}$, and $\Delta U \in \mathbb{U}_{n}$. Then the perturbed
LU factorization of $A$ is
\begin{eqnarray}
A + \Delta A= (L + \Delta L)(U + \Delta U).
\end{eqnarray}
In the following, we regard the perturbations $\Delta L$ and $\Delta U$ as the unknown matrices of the matrix equation (3.1),
and obtain the condition under which the equation (3.1) has the unique solution.

Considering $A=LU$, Eqn. (3.1) can be simplified as
\begin{eqnarray}
L (\Delta U)+ (\Delta L)U=\Delta A- (\Delta L)(\Delta U).
\end{eqnarray}
Premultiplying (3.2) by $L^{-1}$ and postmultiplying it by $U^{-1}$ gives
\begin{eqnarray*}
(\Delta U)U^{-1}+ L^{-1}(\Delta L)=L^{-1}[\Delta A-(\Delta L)(\Delta
U)]U^{-1}.
\end{eqnarray*}
Since $L^{-1}(\Delta L)$ is strictly lower triangular and $(\Delta U)U^{-1}$ is upper triangular, we have
\begin{eqnarray}
&&L^{-1}(\Delta L)={\rm {slt}}\left(L^{-1}[\Delta A-(\Delta L)(\Delta U)]U^{-1}\right),\\
&&(\Delta U)U^{-1}={\rm {ut}}\left(L^{-1}[\Delta A-(\Delta L)(\Delta
U)]U^{-1}\right).
\end{eqnarray}
Let $U_{n-1}$ denote the sub-matrix of $U$ consisting of the first $n-1$ rows and the
first $n-1$ columns, and write $U=\begin{bmatrix}
{U}_{n-1} & u \\
0 & u_{nn} \\
\end{bmatrix}$. Thus, from (3.3), considering the definition of `slt,' it follows that
\begin{eqnarray*}
&&L^{-1}(\Delta L)={\rm {slt}}\left(L^{-1}[\Delta A-(\Delta
L)(\Delta U)]\begin{bmatrix}
{U}_{n-1}^{-1} & 0 \\
0 & 0\\
\end{bmatrix}\right).
\end{eqnarray*}
Applying the operator `vec' to the above equation and using (2.7) and (2.4) implies
\begin{eqnarray*}
(I_n\otimes L^{-1}){\rm {vec}}(\Delta L)=M_{\rm {slt}}\left(\begin{bmatrix}
{U}_{n-1}^{-T} & 0 \\
0 & 0\\
\end{bmatrix}\otimes L^{-1}\right){\rm {vec}}[\Delta A-(\Delta L)(\Delta U)].
\end{eqnarray*}
Premultiplying the above equation by $I_n\otimes L$ and noting (2.9), we get
\begin{eqnarray}
{\rm {vec}}(\Delta L)=(I_n\otimes L)M_{\rm {slt}}\left(\begin{bmatrix}
{U}_{n-1}^{-T} & 0 \\
0 & 0\\
\end{bmatrix}\otimes L^{-1}\right){\rm {vec}}[\Delta A-(\Delta L)(\Delta U)].
\end{eqnarray}
Noticing the structure of $\Delta L$, from (2.4), (2.6), and (2.3), it is seen that
\begin{eqnarray}
{\rm {vec}}(\Delta
L)={\rm {vec}}({\rm {slt}}(\Delta
L))=M_{{\rm {slt}}}{\rm {vec}}(\Delta
L)=M_{{\rm {slvec}}}^TM_{{\rm {slvec}}}{\rm {vec}}(\Delta L)=M_{{\rm
{slvec}}}^T{\rm {slvec}}(\Delta L).
\end{eqnarray}
Substituting the above equality into (3.5) and then left-multiplying it by $M_{{\rm
{slvec}}}$ and using (2.5) yields
\begin{eqnarray}
{\rm {slvec}}(\Delta L)=M_{{\rm
{slvec}}}(I_n\otimes L)M_{\rm {slt}}\left(\begin{bmatrix}
{U}_{n-1}^{-T} & 0 \\
0 & 0\\
\end{bmatrix}\otimes L^{-1}\right){\rm {vec}}[\Delta A-(\Delta L)(\Delta U)].
\end{eqnarray}
Multiplying both sides of (3.7) from the left by $M_{{\rm {slvec}}}^T$ and noting (3.6) and (2.6) leads to
\begin{eqnarray}
{\rm {vec}}(\Delta L)=M_{\rm {slt}}(I_n\otimes L)M_{\rm {slt}}\left(\begin{bmatrix}
{U}_{n-1}^{-T} & 0 \\
0 & 0\\
\end{bmatrix}\otimes L^{-1}\right){\rm {vec}}[\Delta A-(\Delta L)(\Delta
U)].
\end{eqnarray}
From the structure of the matrix $M_{\rm {slt}}$, we can verify that $M_{\rm {slt}}(I_n\otimes L)M_{\rm {slt}}=(I_n\otimes L)M_{\rm {slt}}$, which together with (3.8) gives (3.5). Thus, the equations (3.5) and (3.7) are equivalent.

Similarly, applying the operator `vec' to (3.4) and using (2.7), (2.4), and (2.9), we have
\begin{eqnarray}
{\rm {vec}}(\Delta U)=\left(U^T\otimes I_n\right)M_{\rm
{ut}}\left(U^{-T}\otimes L^{-1}\right){\rm {vec}}[\Delta A-(\Delta
L)(\Delta U)].
\end{eqnarray}
It follows from the structure of $\Delta U$, (2.4), (2.6), and (2.3) that
\begin{eqnarray}
{\rm {vec}}(\Delta
U)={\rm {vec}}({\rm {ut}}(\Delta
U))=M_{{\rm {ut}}}{\rm {vec}}(\Delta
U)=M_{{\rm {uvec}}}^TM_{{\rm {uvec}}}{\rm {vec}}(\Delta U)=M_{{\rm
{uvec}}}^T{\rm {uvec}}(\Delta U).
\end{eqnarray}
Thus, (3.9), (3.10), and (2.5) together implies
\begin{eqnarray}
{\rm {uvec}}(\Delta U)=M_{{\rm {uvec}}}\left(U^T\otimes
I_n\right)M_{\rm {ut}}\left(U^{-T}\otimes L^{-1}\right){\rm
{vec}}[\Delta A-(\Delta L)(\Delta U)].
\end{eqnarray}
Similar to the discussion for $\Delta L$, from (3.11), considering (3.10), (2.6), and the fact $M_{\rm {ut}}(U^T\otimes I_n)M_{\rm {ut}}=(U^T\otimes I_n)M_{\rm {ut}}$, we get (3.9). So the equations (3.9) and (3.11) are equivalent.

Applying the operators `${\rm
{slvec^\dag}}$' and `${\rm {uvec^\dag}}$' to (3.7) and (3.11), respectively, gives
\begin{eqnarray}
&&\Delta L={\rm {slvec^\dag}}\Big(Y_L{\rm {vec}}[\Delta A-(\Delta L)(\Delta U)]\Big),
\end{eqnarray}
and
\begin{eqnarray}
&&\Delta U={\rm {uvec^\dag}}\Big(Y_U{\rm {vec}}[\Delta A-(\Delta L)(\Delta U)]\Big),
\end{eqnarray}
where
\begin{eqnarray*}
Y_L=M_{{\rm
{slvec}}}(I_n\otimes L)M_{\rm {slt}}\left(\begin{bmatrix}
{U}_{n-1}^{-T} & 0 \\
0 & 0\\
\end{bmatrix}\otimes L^{-1}\right),\quad
Y_U=M_{{\rm {uvec}}}\left(U^T\otimes I_n\right)M_{\rm
{ut}}\left(U^{-T}\otimes L^{-1}\right).
\end{eqnarray*}
The matrices $Y_L$ and $Y_U$ are just the ones in \cite[Eqn. (3.5)]{Chang98a}, where their explicit expressions are not given. This fact can be obtained from (3.7) and (3.11),  and \cite[Eqn. (3.6)]{Chang98a} by setting $t=\varepsilon$ in \cite[Eqn. (3.6)]{Chang98a} and dropping the higher-order
terms.

Now we apply the technique of Lyapunov
majorant function and the Banach fixed point principle to derive the rigorous perturbation bounds for $\Delta L$ and $\Delta U$ on the basis of (3.12) and (3.13).

Let $\Delta X=\left[ {\begin{array}{*{20}c}
   \Delta L   \\
   \Delta U  \\
\end{array}} \right]$. Then the equations (3.12) and (3.13) can be rewritten as an operator equation for $\Delta
X$,
\begin{eqnarray}
\Delta X=\Phi(\Delta X,\Delta A)=\left[ {\begin{array}{*{20}c}
   \Phi_1(\Delta X,\Delta A)   \\
   \Phi_2(\Delta X,\Delta A)  \\
\end{array}} \right],
\end{eqnarray}
where
$\Phi_1(\Delta X,\Delta A)={\rm {slvec^\dag}}\Big(Y_L{\rm {vec}}[\Delta A-(\Delta L)(\Delta U)]\Big)$
and
$\Phi_2(\Delta X,\Delta A)={\rm {uvec^\dag}}\Big(Y_U{\rm {vec}}[\Delta A-(\Delta L)(\Delta U)]\Big)$.
Assume that  $Z_1\in {\mathbb{SL}_{n}}$, $Z_2\in {\mathbb{U}_{n}}$, and
$Z=\left[ {\begin{array}{*{20}c}
   Z_1  \\
   Z_2  \\
\end{array}} \right]$. Replacing $\Delta X$ in (3.14) with $Z$ gives
\begin{eqnarray}
Z=\Phi(Z,\Delta A)=\left[ {\begin{array}{*{20}c}
   \Phi_1(Z,\Delta A)   \\
   \Phi_2(Z,\Delta A)  \\
\end{array}} \right],
\end{eqnarray}
where
$\Phi_1(Z,\Delta A)$ and $\Phi_2(Z,\Delta A)$ are the same as $\Phi_1(\Delta X,\Delta A)$ and $\Phi_2(\Delta X,\Delta A)$, respectively, with $\Delta X$ being replaced by $Z$.
Let $|||Z|||\leq \rho=\left[ {\begin{array}{*{20}c}
   \rho_1   \\
   \rho_2  \\
\end{array}} \right]$, i.e, $\left\|
Z_1\right\|_F\leq\rho_1$ and $\left\|
Z_2\right\|_F\leq\rho_2$ for some $\rho_1\geq0$ and $\rho_2\geq0$, and $\left\|\Delta
A\right\|_F= \delta$. Then it follows from the definitions of the
`generalized matrix norm' (2.2) and the operators `${\rm
{uvec^\dag}}$' and `${\rm {slvec^\dag}}$,' with (2.1), that
\begin{eqnarray*}
|||\Phi(Z,\Delta A)|||=\left[ {\begin{array}{*{20}c}
   \left\|\Phi_1(Z,\Delta A)\right\|_F   \\
   \left\|\Phi_2(Z,\Delta A)\right\|_F  \\
\end{array}} \right]\leq\left[ {\begin{array}{*{20}c}
   \left\|Y_L\right\|_2(\delta+\rho_1\rho_2)   \\
   \left\|Y_U\right\|_2(\delta+\rho_1\rho_2) \\
\end{array}} \right].
\end{eqnarray*}
Thus, we have the Lyapunov majorant function (see,
e.g., \cite[Chapter 5]{Konstantinov03}) of
the operator equation (3.15)
\begin{eqnarray*}
h(\rho,\delta)=\left[ {\begin{array}{*{20}c}
   h_1(\rho,\delta)   \\
   h_2(\rho,\delta) \\
\end{array}} \right]=\left[ {\begin{array}{*{20}c}
   \left\|Y_L\right\|_2(\delta+\rho_1\rho_2)   \\
   \left\|Y_U\right\|_2(\delta+\rho_1\rho_2) \\
\end{array}} \right],
\end{eqnarray*}
and the Lyapunov majorant equation (see,
e.g., \cite[Chapter 5]{Konstantinov03})
\begin{eqnarray*}
h(\rho,\delta)=\rho,\ {\textrm{i.e.}, }\quad \left\{
\begin{array}{l}
 \left\|Y_L\right\|_2(\delta+\rho_1\rho_2)=\rho_1, \\
 \left\|Y_U\right\|_2(\delta+\rho_1\rho_2)=\rho_2. \\
 \end{array} \right.
\end{eqnarray*}
Then
\begin{eqnarray}
\rho_2=\frac{\left\|Y_U\right\|_2}{\left\|Y_L\right\|_2}\rho_1,
\end{eqnarray}
and
\begin{eqnarray}
\left\|Y_U\right\|_2\rho_1^2-\rho_1+\left\|Y_L\right\|_2\delta=0.
\end{eqnarray}
Assume that $\delta\in \Omega
=\{\delta\geq0:
1-4\left\|Y_U\right\|_2\left\|Y_L\right\|_2\delta\geq 0\}$.
Then, the Lyapunov majorant equation (3.17) has two nonnegative
roots: $\rho_{1,1}(\delta)\leq\rho_{1,2}(\delta)$ with
\begin{eqnarray*}
\rho_{1,1}(\delta):=f_1(\delta):=\frac{2\left\|Y_L\right\|_
2\delta}{1+\sqrt{1-4\left\|Y_U\right\|_2\left\|Y_L\right\|_2\delta}},
\end{eqnarray*}
which combined with (3.16) gives:
$\rho_{2,1}(\delta)\leq\rho_{2,2}(\delta)$ and
\begin{eqnarray*}
\quad\rho_{2,1}(\delta):=f_2(\delta):=\frac{2\left\|Y_U\right\|_2\delta}{1+\sqrt{1-4\left\|Y_U\right\|_2\left\|Y_L\right\|_2\delta}}.
\end{eqnarray*}
Let the set ${\cal B}(\delta)$ be defined by
\begin{eqnarray*}
{\cal  B}(\delta)=\left\{Z=\left[ {\begin{array}{*{20}c}
   Z_1   \\
   Z_2 \\
\end{array}} \right], Z_1\in {\mathbb{SL}_{n}}, Z_2\in {\mathbb{U}_{n}}: |||Z|||\leq \left[ {\begin{array}{*{20}c}
   f_{1}(\delta)   \\
   f_{2}(\delta) \\
\end{array}} \right]\right\}\subset \mathbb{R}^{2n\times n},
\end{eqnarray*}
which is closed and convex.
Thus, the operator $\Phi(\cdot,\Delta A)$ maps the set ${\cal
B}(\delta)$ into itself. Furthermore, note that the Jacobi matrix of
$h(\rho,\delta)$ relative to $\rho$ at $\rho_0$ is,
\begin{eqnarray*}
h^{'}_\rho(\rho_0,\delta)=\frac{1-\sqrt{1-4\left\|Y_U\right\|_2\left\|Y_L\right\|_2\delta}}{2}\left[ {\begin{array}{*{20}c}
   1 & \left\|Y_L\right\|_2/\left\|Y_U\right\|_2  \\
   \left\|Y_U\right\|_2/\left\|Y_L\right\|_2 &   1\\
\end{array}} \right],
\end{eqnarray*}
where $\rho_0=\left[ {\begin{array}{*{20}c}
   f_1(\delta)   \\
   f_2(\delta) \\
\end{array}} \right]$, and for $Z, \widetilde{Z}\in {\cal B}(\delta)$,
\begin{eqnarray*}
||| \Phi(Z,\Delta A)-\Phi(\widetilde{Z},\Delta A) ||| \leq
h^{'}_\rho(\rho_0,\delta)|||Z-\widetilde{Z}|||.
\end{eqnarray*}
Then if $\delta\in \Omega_1
=\{\delta\geq0: 1-4\left\|Y_U\right\|_2\left\|Y_L\right\|_2\delta>
0\}$,
we have that the spectral radius of $h^{'}_\rho(\rho_0,\delta)$ is
smaller than 1 and then the operator $\Phi(\cdot,\Delta A)$ is
generalized contractive (see,
e.g., \cite[Appendix D]{Konstantinov03}) on ${\cal B}(\delta)$.
According to the generalized Banach fixed point principle (see,
e.g., \cite[Appendix D]{Konstantinov03}), there exists a unique solution to the
operator equation (3.15) in the set ${\cal B}(\delta)$ when $\delta\in \Omega_1$, and so does the operator equation (3.14). As a result,
we have
\begin{eqnarray*}
|||\Delta X|||\leq\left[ {\begin{array}{*{20}c}
   f_1(\delta)   \\
   f_2(\delta) \\
\end{array}} \right],\quad \delta \in \Omega_1.
\end{eqnarray*}
Considering the equivalence of the matrix equation (3.1) and the
operator equation (3.14),  we have the main theorem.

\begin{theorem}
Let the unique LU factorization of $A\in \mathbb{R}^{n\times n}_n$
be as in {\rm (1.1)} and  $\Delta A\in \mathbb{R}^{n\times n}$. If
\begin{eqnarray}
\left\|Y_L\right\|_2\left\|Y_U\right\|_2\left\|\Delta A\right\|_F<\frac{1}{4},
\end{eqnarray}
then $A+\Delta A$ has the unique LU factorization {\rm (3.1)}. Moreover,
\begin{eqnarray}
&&\left\|\Delta L\right\|_F\leq\frac{2\left\|Y_L\right\|_2\left\|\Delta A\right\|_F}{1+\sqrt{1-4\left\|Y_U\right\|_2\left\|Y_L\right\|_2\left\|\Delta A\right\|_F}}\\
&&\quad\quad\quad \leq 2\left\|Y_L\right\|_2\left\|\Delta
A\right\|_F\\
&&\quad\quad\quad = 2\left\|(I_n\otimes L)M_{\rm {slt}}\left(\begin{bmatrix}
{U}_{n-1}^{-T} & 0 \\
0 & 0\\
\end{bmatrix}\otimes L^{-1}\right)\right\|_2\left\|\Delta
A\right\|_F,
\end{eqnarray}
and
\begin{eqnarray}
&&\left\|\Delta U\right\|_F\leq\frac{2\left\|Y_U\right\|_2\left\|\Delta A\right\|_F}{1+\sqrt{1-4\left\|Y_U\right\|_2\left\|Y_L\right\|_2\left\|\Delta A\right\|_F}}\\
&&\quad\quad\quad \leq 2\left\|Y_U\right\|_2\left\|\Delta A\right\|_F\\
&&\quad\quad\quad = 2\left\|\left( U^T\otimes
I_n\right)M_{\rm {ut}}\left(U^{-T}\otimes
L^{-1}\right)\right\|_2\left\|\Delta A\right\|_F.
\end{eqnarray}
\end{theorem}
{\raggedleft\em  Proof.}
From the discussions before Theorem 3.1, we only need to show that (3.21) and (3.24) hold.
Considering the definition of the spectral norm, (2.6), and the facts
\begin{eqnarray}
M_{\rm {slt}}(I_n\otimes L)M_{\rm {slt}}=(I_n\otimes L)M_{\rm {slt}},\quad M_{\rm {ut}}(U^T\otimes I_n)M_{\rm {ut}}=(U^T\otimes I_n)M_{\rm {ut}},
\end{eqnarray}
it is easy to verify that
\begin{eqnarray*}
\left\|Y_L\right\|_2=\left\|(I_n\otimes L)M_{\rm {slt}}\left(\begin{bmatrix}
{U}_{n-1}^{-T} & 0 \\
0 & 0\\
\end{bmatrix}\otimes L^{-1}\right)\right\|_2,\quad \left\|Y_U\right\|_2=\left\|\left( U^T\otimes
I_n\right)M_{\rm {ut}}\left(U^{-T}\otimes
L^{-1}\right)\right\|_2.
\end{eqnarray*}
So (3.21) and (3.24) hold.
 $\square$

\begin{remark} {\em   From (3.19) and (3.22), we have the following first-order
perturbation bounds,
\begin{eqnarray}
&&\left\|\Delta L\right\|_F\leq \left\|(I_n\otimes L)M_{\rm {slt}}\left(\begin{bmatrix}
{U}_{n-1}^{-T} & 0 \\
0 & 0\\
\end{bmatrix}\otimes L^{-1}\right)\right\|_2\left\|\Delta
A\right\|_F+ {\cal O}\left(\left\|\Delta A\right\|_F^2\right),
\end{eqnarray}
and
\begin{eqnarray}
&&\left\|\Delta U\right\|_F\leq \left\|\left( U^T\otimes
I_n\right)M_{\rm {ut}}\left(U^{-T}\otimes
L^{-1}\right)\right\|_2\left\|\Delta A\right\|_F+{\cal O}\left(\left\|\Delta A\right\|_F^2\right).
\end{eqnarray}
Note that, in this case, the condition (3.18) can be weakened to
\begin{eqnarray}
\left\|L^{-1}\right\|_2\left\|U^{-1}\right\|_2\left\|\Delta A\right\|_F<1.
\end{eqnarray}
This is because the bounds (3.26) and (3.27) can be derived from
(3.12) and (3.13) directly by omitting the higher-order
terms. We only provide the condition under which the LU
factorization of $A+\Delta A$ exists and is unique. From \cite[Proof of Theorem 4.1]{Chang10},
it follows that the condition (3.28) is enough.

The bounds (3.26) and (3.27) without explicit expressions were also derived by the matrix-vector equation
approach in \cite{Chang98a}, which are considered to be optimal.}
\end{remark}

\begin{remark} {\em   The rigorous perturbation bounds derived by the
combination of the classic and refined matrix equation approaches
presented in \cite{Chang10} are as follows,
\begin{eqnarray}
&&\left\|\Delta L\right\|_F\leq 2\left(\mathop {\inf }\limits_{D_L
\in {\mathbb{D}_{n}}} {{k_2}\left({{LD_L^{ - 1}}} \right)}
\right)\left\|U^{-1}_{n-1}\right\|_2\left\|\Delta A\right\|_F,
\end{eqnarray}
and
\begin{eqnarray}
&&\left\|\Delta U\right\|_F\leq 2\left(\mathop
{\inf }\limits_{D_U \in {\mathbb{D}_{n}}} {{k_2}\left({{D_U^{ -
1}U}} \right)} \right)\left\|L^{-1}\right\|_2\left\|\Delta
A\right\|_F,
\end{eqnarray}
under the condition
\begin{eqnarray}
\left\|L^{-1}\right\|_2\left\|U^{-1}\right\|_2\left\|\Delta A\right\|_F< 1/4.
\end{eqnarray}
In (3.29) and (3.30), $\mathbb{D}_{n}$ denotes the set of $n
\times n$ positive definite diagonal matrices. The bounds (3.29)
and (3.30) can be much smaller than the previous ones derived by the
classic matrix equation approach; see discussions in
\cite{Chang10}. From \cite[Eqns. (3.17) and (3.24)]{Chang98a}, we
have
\begin{eqnarray*}
\left\|Y_L\right\|_2\leq \left(\mathop {\inf }\limits_{D_L \in {\mathbb{D}_{n}}} {{k_2}\left( {{LD_L^{ - 1}}} \right)}\right)\left\|U^{-1}_{n-1}\right\|_2,\quad
\left\|Y_U\right\|_2\leq \left(\mathop {\inf }\limits_{D_U \in {\mathbb{D}_{n}}} {{k_2}\left( {{D_U^{ - 1}U}} \right)} \right)\left\|L^{-1}\right\|_2.
\end{eqnarray*}
So the bounds (3.21) and (3.24) are tighter than (3.29) and (3.30),
respectively. Unfortunately, it follows from \cite[Eqns.
(3.18) and (3.25)]{Chang98a} that
\begin{eqnarray*}
\left\|Y_L\right\|_2\geq \left\|U^{-1}_{n-1}\right\|_2,\quad\left\|Y_U\right\|_2\geq \left\|L^{-1}\right\|_2.
\end{eqnarray*}
Thus, the condition (3.18) is more constraining than (3.31).
Fortunately, the above two lower bounds are attainable
\cite{Chang97a, Chang98a}, which shows that the condition
(3.18) is not so constraining. In addition, it is also a little more
expensive to estimate the bounds (3.21) and (3.24) than that of
(3.29) and (3.30) because the former involves the Kronecker products.
These should be the price of having tighter rigorous perturbation
results.}
\end{remark}

Considering the standard techniques of backward error analysis
(see, e.g., \cite[Theorem 9.3]{Higham02}), we have that the computed LU factors
$\widetilde{L}$ and $\widetilde{U}$ by the Gaussian elimination
satisfy,
\begin{eqnarray}
\widetilde{A}=A+\Delta A=\widetilde{L}\widetilde{U},\quad |\Delta
A|\leq \varepsilon
|\widetilde{L}||\widetilde{U}|,
\end{eqnarray}
where $\varepsilon=n {\bf u}/(1-n {\bf u})$ with ${\bf u}$ being the
unit roundoff. In the following, we consider the rigorous
perturbation bounds for the LU factorization with the perturbation
$\Delta A$ having the same form as in (3.32). The new bounds, similar to
the ones in \cite{Chang10}, will involve the LU factors of
$\widetilde{A}$. The reader can refer to \cite[Section 4]{Chang10} for an
explanation.

Assume that the matrices $\widetilde{A}$, $\widetilde{L}$, and $\widetilde{U}$ in (3.32) are perturbed as
\begin{eqnarray*}
\widetilde{A}\rightarrow \widetilde{A}-\Delta A,\ \widetilde{L}\rightarrow \widetilde{L}-\Delta L,\ \widetilde{U}\rightarrow \widetilde{U}-\Delta U,
\end{eqnarray*}
where $\Delta A\in {\mathbb{R}^{n \times n}}$ is as in (3.32),
$\Delta L\in {\mathbb{SL}_{n}}$, and $\Delta U \in {\mathbb{U}_{n
}}$. Then the perturbed LU factorization of $\widetilde{A}$ is
\begin{eqnarray*}
A=\widetilde{A}-\Delta A=(\widetilde{L}-\Delta
L)(\widetilde{U}-\Delta U),
\end{eqnarray*}
which together with (3.32) yields,
\begin{eqnarray*}
\widetilde{L}(\Delta U) + (\Delta L)\widetilde{U}=\Delta A+ (\Delta
L)(\Delta U).
\end{eqnarray*}
As done before, we regard the perturbations $\Delta{L}$ and $\Delta{U}$ as the unknown matrices. Thus, similar to the induction before Theorem 3.1,
replacing $L$ and $U$ with $\widetilde{L}$ and $\widetilde{U}$, respectively, we have
\begin{eqnarray}
\Delta X=\widetilde{\Phi}(\Delta X,\Delta A)=\left[ {\begin{array}{*{20}c}
   \widetilde{\Phi}_1(\Delta X,\Delta A)   \\
   \widetilde{\Phi}_2(\Delta X,\Delta A)  \\
\end{array}} \right],
\end{eqnarray}
where
\begin{eqnarray}
\widetilde{\Phi}_1(\Delta X,\Delta A)={\rm
{slvec^\dag}}\Big(Y_{\widetilde{L}}{\rm {vec}}(\Delta
A)+Y_{\widetilde{L}}{\rm {vec}}[(\Delta L)(\Delta U)]\Big)
\end{eqnarray}
and
\begin{eqnarray}
\widetilde{\Phi}_2(\Delta X,\Delta A)={\rm
{uvec^\dag}}\Big(Y_{\widetilde{U}}{\rm {vec}}(\Delta
A)+Y_{\widetilde{U}}{\rm {vec}}[(\Delta L)(\Delta U)]\Big).
\end{eqnarray}
Here
\begin{eqnarray*}
Y_{\widetilde{L}}=M_{{\rm {slvec}}}\left(I_n\otimes
\widetilde{L}\right)M_{\rm {slt}}\left(\begin{bmatrix}
{\widetilde{U}}_{n-1}^{-T} & 0 \\
0 & 0\\
\end{bmatrix}\otimes \widetilde{L}^{-1}\right),\quad
Y_{\widetilde{U}}=M_{{\rm {uvec}}}\left(\widetilde{U}^T\otimes
I_n\right)M_{\rm {ut}}\left(\widetilde{U}^{-T}\otimes
\widetilde{L}^{-1}\right).
\end{eqnarray*}

Considering (3.32), the fact that the Frobenius norm is monotone,
and (2.1), we obtain
\begin{eqnarray*}
\left\|\widetilde{\Phi}_1(Z,\Delta A)\right\|_F\leq\left\|
|Y_{\widetilde{L}}|{\rm
{vec}}\left(|\widetilde{L}|
|\widetilde{U}|\right)\right\|_F\varepsilon+\left\|
|Y_{\widetilde{L}}|\right\|_2\rho_1\rho_2
\end{eqnarray*}
and
\begin{eqnarray*}
\left\|\widetilde{\Phi}_2(Z,\Delta A)\right\|_F\leq\left\|
|Y_{\widetilde{U}}|{\rm
{vec}}\left(|\widetilde{L}|
|\widetilde{U}|\right)\right\|_F\varepsilon+\left\|
|Y_{\widetilde{U}}|\right\|_2\rho_1\rho_2.
\end{eqnarray*}
Similar to the discussions before Theorem 3.1, using the above two inequalities, we have the following theorem.

\begin{theorem}
Assume that $\Delta A\in \mathbb{R}^{n\times n}$ is a perturbation in $A\in \mathbb{R}^{n\times n}$ and $A+\Delta A$ has the unique LU factorization satisfying {\rm (3.32)}. Let $a=\left\||Y_{\widetilde{L}}|{\rm {vec}}(|\widetilde{L}||\widetilde{U}|)\right\|_F$, $b=\left\||Y_{\widetilde{U}}|{\rm {vec}}(|\widetilde{L}||\widetilde{U}|)\right\|_F$, and $c=b\left\| |Y_{\widetilde{L}}|\right\|_2-a\left\|
|Y_{\widetilde{U}}|\right\|_2$. If
\begin{eqnarray}
|c|\varepsilon<1\ \textrm{ and }\
4a\left\||Y_{\widetilde{U}}|\right\|_2\varepsilon<\left(1-c\varepsilon\right)^2,
\end{eqnarray}
then $A$ has the unique LU factorization $A=LU$, where $L=\widetilde{L} -\Delta L$ and $U=\widetilde{U}-\Delta U$. Moreover,
\begin{eqnarray}
&&\left\|\Delta L\right\|_F\leq\frac{2a\varepsilon}{1-c
\varepsilon+\sqrt{\left(1-c
\varepsilon\right)^2-4a\left\||Y_{\widetilde{U}}|\right\|_2\varepsilon}}\\
&&\quad\quad\quad\
\leq\frac{2a\varepsilon}{1-c\varepsilon},
\end{eqnarray}
and
\begin{eqnarray}
&&\left\|\Delta U\right\|_F\leq\frac{2b\varepsilon}{1+c
\varepsilon+\sqrt{\left(1-c\varepsilon\right)^2
-4a\left\||Y_{\widetilde{U}}|\right\|_2\varepsilon}}\\
&&\quad\quad\quad\ \leq\frac{2b\varepsilon}{1+c\varepsilon}.
\end{eqnarray}
\end{theorem}
\begin{remark} {\em   From (3.37) and (3.39), we have the following
first-order perturbation bounds,
\begin{eqnarray}
\left\|\Delta L\right\|_F\leq
\left\||Y_{\widetilde{L}}|{\rm
{vec}}\left(|\widetilde{L}||\widetilde{U}|\right)\right\|_F\varepsilon+{\cal
O}\left(\varepsilon^2 \right),
\end{eqnarray}
and
\begin{eqnarray}
\left\|\Delta
U\right\|_F\leq \left\||Y_{\widetilde{U}}|{\rm
{vec}}\left(|\widetilde{L}|
|\widetilde{U}|\right)\right\|_F\varepsilon+{\cal
O}\left(\varepsilon^2\right),
\end{eqnarray}
which can also be derived from (3.33)--(3.35), and (3.32) directly
by omitting the higher-order terms. Therefore, in this case, the
condition (3.36) can be weakened to
\begin{eqnarray}
\left\||\widetilde{L}^{-1}||\widetilde{L}|\right\|_F\left\||\widetilde{U}|
|\widetilde{U}^{-1}|\right\|_F\varepsilon< 1,
\end{eqnarray}
which guarantees that the unique LU factorization of
$\widetilde{A}-\Delta A=A$ exists \cite[Proof of Theorem
4.2]{Chang10}. Using (3.33)--(3.35), and (3.32), we can also obtain
the first-order perturbation bounds with respect to the `$M$'-norm
and the `$S$'-norm,
\begin{eqnarray}
\left\|\Delta L\right\|_{\nu}\leq
\left\||Y_{\widetilde{L}}|{\rm
{vec}}\left(|\widetilde{L}||\widetilde{U}|\right)\right\|_{\nu}\varepsilon+{\cal
O}\left(\varepsilon^2\right),
\end{eqnarray}
and
\begin{eqnarray}
\left\|\Delta U\right\|_{\nu}\leq
\left\||Y_{\widetilde{U}}|{\rm
{vec}}\left(|\widetilde{L}||\widetilde{U}|\right)\right\|_{\nu}\varepsilon+{\cal
O}\left(\varepsilon^2\right),
\end{eqnarray}
where $\nu=M$ or $S$, under the condition (3.43). Recall that the
$M$-norm and the $S$-norm of a matrix $A=(a_{ij})\in
\mathbb{R}^{m\times n}$ are defined by (see, e.g., \cite[Chapter 6]{Higham02}),
\begin{eqnarray*}
\left\|A\right\|_M =\mathop {\max }\limits_{i,j}|a_{ij}|,\quad\left\|A\right\|_S =\mathop {\sum}\limits_{i,j}|a_{ij}|,
\end{eqnarray*}
respectively, which are both monotone. For the $M$-norm, the
first-order bound for $\Delta L$, i.e., (3.44), is attained for $\Delta A$
satisfying
\begin{eqnarray*}
{\rm vec}(\Delta A)=\varepsilon D_k{\rm
vec}\left(|\widetilde{L}|
|\widetilde{U}|\right),\quad D_k={\rm
diag}\left(\xi_1,\xi_2,\cdots,\xi_{n^2}\right),
\end{eqnarray*}
where $\xi_i={\rm sign}\left(Y_{\widetilde{L}}(k,i)\right)$ and
$\left\||Y_{\widetilde{L}}|{\rm
{vec}}\left(|\widetilde{L}||\widetilde{U}|\right)\right\|_M=\left(|Y_{\widetilde{L}}|{\rm
{vec}}\left(|\widetilde{L}||\widetilde{U}|\right)\right)(k,1)$.
Here, the MATLAB notation is used. If we take $\xi_i={\rm
sign}\left(Y_{\widetilde{U}}(k,i)\right)$ and
$\left\||Y_{\widetilde{U}}|{\rm
{vec}}\left(|\widetilde{L}||\widetilde{U}|\right)\right\|_M=\left(|Y_{\widetilde{U}}|{\rm
{vec}}\left(|\widetilde{L}||\widetilde{U}|\right)\right)$
$(k,1)$,
then the first-order bound for $\Delta U$, i.e., (3.45), is attained under
the $M$-norm for this $\Delta A$. Thus, we obtain the optimal
first-order perturbation bounds for the LU factorization under the
$M$-norm.

In \cite{Chang02}, Chang presented the following first-order perturbation
bounds under the consistent and monotone norm
$\left\|\cdot\right\|$,
\begin{eqnarray}
\left\|\Delta L\right\|\leq
\left\||\widetilde{L}||\widetilde{L}^{-1}||\widetilde{L}|\right\|\cdot\left\|
|\widetilde{U}_{n-1}||\widetilde{U}^{-1}_{n-1}|\right\|\varepsilon+{\cal
O}\left(\varepsilon^2 \right),
\end{eqnarray}
and
\begin{eqnarray}
\left\|\Delta U\right\|\leq
\left\||\widetilde{U}||\widetilde{U}^{-1}||\widetilde{U}|\right\|\cdot\left\|
|\widetilde{L}^{-1}||\widetilde{L}|\right\|\varepsilon+{\cal
O}\left(\varepsilon^2\right).
\end{eqnarray}
Since, for the norm $\left\|\cdot\right\|_{\nu}$ ($\nu=F$ or $S$), which are both consistent
and monotone, considering (2.7), (2.4), and (2.3), we have
\begin{eqnarray}
&&\left\||Y_{\widetilde{L}}|{\rm
{vec}}\left(|\widetilde{L}||\widetilde{U}|\right)\right\|_{\nu}\leq
\left\|M_{{\rm {slvec}}}\left(I_n\otimes
|\widetilde{L}|\right)M_{\rm {slt}}\left(\begin{bmatrix}
|{\widetilde{U}}_{n-1}^{-T}| & 0 \\
0 & 0\\
\end{bmatrix}\otimes |\widetilde{L}^{-1}|\right){\rm
{vec}}\left(|\widetilde{L}||\widetilde{U}|\right)\right\|_{\nu}\nonumber\\
&&\quad\quad\quad\quad\quad\quad\quad\quad\ \ \ =\left\|{\rm
{slvec}}\left(|\widetilde{L}|{\rm
{slt}}\left(|\widetilde{L}^{-1}||\widetilde{L}||\widetilde{U}|\begin{bmatrix}
|{\widetilde{U}}_{n-1}^{-1}| & 0 \\
0 & 0\\
\end{bmatrix}\right)\right)\right\|_{\nu}\nonumber\\
&&\quad\quad\quad\quad\quad\quad\quad\quad\ \ \
=\left\||\widetilde{L}|{\rm
{slt}}\left(|\widetilde{L}^{-1}||\widetilde{L}||\widetilde{U}|\begin{bmatrix}
|\widetilde{U}_{n-1}^{-1}| & 0 \\
0 & 0 \\
\end{bmatrix}\right)\right\|_{\nu}\\
&&\quad\quad\quad\quad\quad\quad\quad\quad\ \ \
\leq\left\||\widetilde{L}||\widetilde{L}^{-1}||\widetilde{L}|\right\|_{\nu}\left\||\widetilde{U}_{n-1}|
|\widetilde{U}^{-1}_{n-1}|\right\|_{\nu},\nonumber
\end{eqnarray}
and
\begin{eqnarray}
&&\left\||Y_{\widetilde{U}}|{\rm {vec}}\left(|\widetilde{L}||\widetilde{U}|\right)\right\|_{\nu}\leq\left\|M_{{\rm {uvec}}}
\left(|\widetilde{U}^T|\otimes I_n\right)M_{\rm {ut}}\left(|\widetilde{U}^{-T}|\otimes |\widetilde{L}^{-1}|\right)
{\rm {vec}}\left(|\widetilde{L}||\widetilde{U}|\right)\right\|_{\nu}\nonumber\\
&&\quad\quad\quad\quad\quad\quad\quad\quad\quad =\left\|{\rm {uvec}}\left({\rm {ut}}\left(|\widetilde{L}^{-1}||\widetilde{L}|
|\widetilde{U}||\widetilde{U}^{-1}|\right)|\widetilde{U}|\right)\right\|_{\nu}\nonumber\\
&&\quad\quad\quad\quad\quad\quad\quad\quad\quad = \left\|{\rm
{ut}}\left(|\widetilde{L}^{-1}||\widetilde{L}||\widetilde{U}||
\widetilde{U}^{-1}|\right)|\widetilde{U}|\right\|_{\nu}\\
&&\quad\quad\quad\quad\quad\quad\quad\quad\quad
\leq\left\||\widetilde{U}||\widetilde{U}^{-1}||\widetilde{U}|\right\|_{\nu}
\left\||\widetilde{L}^{-1}||\widetilde{L}|\right\|_{\nu},\nonumber
\end{eqnarray}
the first-order bounds (3.44) and (3.45) are tighter than (3.46)
and (3.47) under the two norms, respectively.

In addition, it should be pointed out that we can not achieve the
first-order perturbation bounds in terms of the 1-norm and the
$\infty$-norm, both of which are also consistent and monotone.}
\end{remark}

\begin{remark} {\em   In \cite{Chang10}, the following rigorous perturbation bounds
with respect to the consistent and monotone norm were derived by the
combination of the classic and refined matrix equation approaches,
\begin{eqnarray}
\left\|\Delta L\right\|\leq 2\left(\mathop {\inf }\limits_{D_L \in
{\mathbb{D}_{n}}} {\left\|{{\widetilde{L}D_L^{ - 1}}} \right\| \cdot
\left\|D_L|\widetilde{L}^{-1}||\widetilde{L}|
\right\|}
 \right)\left\||\widetilde{U}_{n-1}| \cdot |\widetilde{U}^{-1}_{n-1}|\right\|\varepsilon,
\end{eqnarray}
and
\begin{eqnarray}
\left\|\Delta U\right\|\leq 2\left(\mathop {\inf }\limits_{D_U \in
{\mathbb{D}_{n}}} {\left\|{{D_U^{ - 1}\widetilde{U}}} \right
\|\cdot\left\||\widetilde{U}||\widetilde{U}^{-1}|D_U\right\|}
\right)\left\||\widetilde{L}^{-1}|\cdot
|\widetilde{L}|\right\|\varepsilon,
\end{eqnarray}
under the condition
\begin{eqnarray}
\left\||\widetilde{L}^{-1}|
|\widetilde{L}|\right\| \cdot
\left\||\widetilde{U}||\widetilde{U}^{-1}|\right
\|\varepsilon< 1/4.
\end{eqnarray}
Combining the properties of the operators `ut' and `slt'
\cite[Eqn (2.5)]{Chang10}
\begin{eqnarray*}
{\rm {slt}}(D_LX)=D_L{\rm {slt}}(X),\quad {\rm {ut}}(XD_U)={\rm {ut}}(X)D_U,
\end{eqnarray*}
where $D_L,D_U\in\mathbb{D}_{n}$, with (3.48) and (3.49), and noting (2.1), we have
\begin{eqnarray*}
&&\left\||Y_{\widetilde{L}}|{\rm
{vec}}\left(|\widetilde{L}||\widetilde{U}|\right)\right\|_F\leq\left\||\widetilde{L}|D_L^{-1}{\rm
{slt}}\left(D_L|\widetilde{L}^{-1}||\widetilde{L}||\widetilde{U}|\begin{bmatrix}
|\widetilde{U}_{n-1}^{-1}| & 0 \\
0 & 0 \\
\end{bmatrix}\right)\right\|_F\nonumber\\
&&\quad\quad\quad\quad\quad\quad\quad\quad\ \ \ \leq\left\|{{|\widetilde{L}|D_L^{ - 1}}}
\right\|_2\left\|D_L|\widetilde{L}^{-1}||\widetilde{L}|
\right\|_2\left\||\widetilde{U}_{n-1}||\widetilde{U}^{-1}_{n-1}|\right\|_F,
\end{eqnarray*}
and
\begin{eqnarray*}
&&\left\||Y_{\widetilde{U}}|{\rm {vec}}\left(|\widetilde{L}||\widetilde{U}|\right)\right\|_F\leq
\left\|{\rm {ut}}\left(|\widetilde{L}^{-1}||\widetilde{L}||\widetilde{U}|| \widetilde{U}^{-1}|D_U\right)
D_U^{-1}|\widetilde{U}|\right\|_F\nonumber\\
&&\quad\quad\quad\quad\quad\quad\quad\quad\ \ \ \ \leq\left\|{{D_U^{ -
1}|\widetilde{U}|}} \right\|_2
\left\||\widetilde{U}||\widetilde{U}^{-1}|D_U\right\|_2\left\||\widetilde{L}^{-1}||\widetilde{L}|\right\|_F.
\end{eqnarray*}
Note that $D_L,D_U\in\mathbb{D}_{n}$ are arbitrary. Thus, under the Frobenius norm, when
\begin{eqnarray*}
\left\||\widetilde{L}|D_L^{ - 1}
\right\|_2=\left\|\widetilde{L}D_L^{ - 1}\right\|_2,\quad
\left\|D_U^{ - 1}|\widetilde{U}| \right\|_2=\left\|D_U^{
- 1}\widetilde{U}\right\|_2,
\end{eqnarray*}
if $-1<c\varepsilon<0$, the bound (3.38) is obviously
smaller than (3.50); if
$1>c\varepsilon>0$,
the bound (3.40) is obviously smaller than (3.51); otherwise, the
bounds (3.38) and (3.40) are obviously smaller than the
corresponding ones (3.50) and (3.51). Notice that for any matrix $X
\in {\mathbb{R}^{m \times n}}$, $\left\||X|\right\|_2$ is at most
$\sqrt{{\rm rank}(X)}$ times as large as $\left\|X\right\|_2$
(see, e.g., \cite[Lemma 6.6]{Higham02}). Especially, the scaling
matrices can make $\widetilde{L}D_L^{ - 1}$ and $D_U^{ -
1}\widetilde{U}$ be of special structure. For example, they may have
the unit 2-norm columns and rows, respectively. As a result, the
differences between $\left\||\widetilde{L}|D_L^{ - 1}
\right\|_2$ and $\left\|\widetilde{L}D_L^{ - 1}\right\|_2$,
$\left\|D_U^{ - 1}|\widetilde{U}| \right\|_2$ and
$\left\|D_U^{ - 1}\widetilde{U}\right\|_2$ will not be remarkable in
general. See the following example. Moreover, since $\varepsilon$ is
very small,
$c\varepsilon$ may also be very small. See Example 3.1 below. Thus, the bounds
(3.38) and (3.40) may generally be smaller than (3.50) and (3.51),
respectively. An example is given below to indicate this conjecture.
However, it should be mentioned that the condition (3.36) is more
complicated and may be more constraining than the one (3.52), and it
is a slightly more expensive to estimate the bounds in Theorem 3.2.

In addition, we need to point out that we can not obtain the
rigorous perturbation bounds under the $S$-norm, the $M$-norm, the
1-norm, and the $\infty$-norm using the foregoing approach.}
\end{remark}

\begin{example} {\em   The example is from \cite{Chang98a}. That is, each test
matrix has the form $A=D_1BD_2$, where $D_1={\rm
diag}(1,d_1,d_1^2,\cdots,d_1^{n-1})$, $D_2={\rm
diag}(1,d_2,d_2^2,\cdots,d_2^{n-1})$, and $B \in {\mathbb{R}^{n
\times n}}$ is a random matrix produced by the MATLAB function $\mathbf{randn}$. As done
in \cite{Chang98a}, the chosen scaling matrices $D_L$ and $D_U$ are defined by
$D_L={\rm diag}(\left\|L(:,j)\right\|_2)$ and $D_U={\rm
diag}(\left\|U(j,:)\right\|_2)$, respectively. Upon
computations in MATLAB 7.0 on a PC, with machine precision $2.2\times 10^{-16}$, the numerical results for $n=10$, $d_1,d_2\in
\{0.2,1,2\}$, and the same matrix $B$ are listed in Table 1, which demonstrate the
conjectures given in Remark 3.4.

\begin{center}
 \begin{tabular}{|r|r|r|r|r|r|r|r|r|r|r|}
    \multicolumn{11}{c}{\centering  Table 1: Comparison of rigorous bounds for $A=D_1BD_2$ } \\\hline
$d_1$& $d_2$ & $\gamma_{_L}$ \quad\  &$\gamma_{_L}(D_L)$ & $\eta_{_{D_L}}$ &$\gamma_{_U} $ \quad\ &$\gamma_{_U}(D_U) $&$\eta_{_{D_U} }$ &$t_{_{\gamma}} $\quad\quad &$t_{_{\gamma(D)}} $ &$\tau $\quad\quad\ \\\hline
0.2  &0.2  &4.31{\rm e}+01  &2.66{\rm e}+06  &1.00  &1.00{\rm e}+00  &5.93{\rm e}+00  &1.00 &0.007 &0.002 &9.15{\rm e}-05 \\
0.2  &1    &4.31{\rm e}+01  &2.66{\rm e}+06  &1.00  &1.38{\rm e}+00  &2.83{\rm e}+02  &1.20 &0.009 &0.001 &6.99{\rm e}-09 \\
0.2  &2    &4.31{\rm e}+01  &2.66{\rm e}+06  &1.00  &1.49{\rm e}+00  & 9.23{\rm e}+02 &1.09 &0.026 &0.002 &5.81{\rm e}-07 \\
1    &0.2  &7.17{\rm e}+01  &6.17{\rm e}+02  &1.27  &1.03{\rm e}+00  &9.13{\rm e}+01  &1.00 &0.010 &0.002 &2.58{\rm e}-09  \\
1    &1    &7.17{\rm e}+01  &6.17{\rm e}+02  &1.27  &1.72{\rm e}+02  &1.68{\rm e}+03  &1.20 &0.018 &0.002 &9.45{\rm e}-11  \\
1    &2    &7.17{\rm e}+01  &6.17{\rm e}+02  &1.27  &2.27{\rm e}+02  &2.65{\rm e}+03  &1.09 &0.014 &0.002 &4.85{\rm e}-08  \\
2    &0.2  &1.27{\rm e}+01  &1.04{\rm e}+03  &1.11  &3.11{\rm e}+00  &1.52{\rm e}+04  &1.00 &0.009 &0.002 &-3.21{\rm e}-09  \\
2    &1    &1.27{\rm e}+01  &1.04{\rm e}+03  &1.11  & 2.79{\rm e}+02  &3.05{\rm e}+04 &1.20 &0.021 &0.003 &1.17{\rm e}-06 \\
2    &2    &1.27{\rm e}+01  &1.04{\rm e}+03  &1.11  &2.78{\rm e}+02  &3.84{\rm e}+04  &1.09 &0.016 &0.002 &3.75{\rm e}-04\\
\hline
\end{tabular}
\end{center}

In Table 1, we denote
\begin{align*}
&\gamma_{_L}=\frac{a}{1-c\varepsilon}/
\left\|\widetilde{L}\right\|_F,\quad \gamma_{_L}(D_L)=\left({\left\|{{\widetilde{L}D_L^{ - 1}}} \right\|_2\left\|D_L|\widetilde{L}^{-1}||\widetilde{L}| \right\|_2} \right)
\left\||\widetilde{U}_{n-1}||\widetilde{U}^{-1}_{n-1}|\right\|_F/\left\|\widetilde{L}\right\|_F,\\
&\gamma_{_U}=\frac{b}{1+c
\varepsilon}/\left\|\widetilde{U}\right\|_F,\quad \gamma_{_U}(D_U)=\left({\left\|{{D_U^{ - 1}\widetilde{U}}} \right\|_2\left\||\widetilde{U}||\widetilde{U}^{-1}|D_U\right\|_2} \right)
\left\||\widetilde{L}^{-1}||\widetilde{L}|\right\|_F/\left\|\widetilde{U}\right\|_F,\\
&\eta_{_{D_L}}=\left\||\widetilde{L}|D_L^{ -
1}\right\|_2/\left\|\widetilde{L}D_L^{ - 1}\right\|_2,\quad
\eta_{_{D_U}}=\left\|D_U^{ -
1}|\widetilde{U}|\right\|_2/\left\|D_U^{ -
1}\widetilde{U}\right\|_2,\quad
\tau=c\varepsilon,
\end{align*}
and $t_{_{\gamma}}$ and $t_{_{\gamma(D)}}$ the time cost for computing $\gamma_{_L},\gamma_{_U}$ and $\gamma_{_L}(D_L),\gamma_{_U}(D_U)$, respectively.
}
\end{example}

\begin{remark}{\em
Considering the definitions of the matrix norms used above, the fact that for any matrix $X\in {\mathbb{R}^{n^2 \times n^2 }}$,
$
|M_{{\rm {uvec}}}X|=M_{{\rm {uvec}}}|X|,\ | M_{{\rm
{slvec}}}X|= M_{{\rm
{slvec}}}|X|,\ |M_{\rm {ut}}X|=M_{\rm {ut}}|X|, \ |M_{\rm {slt}}X|=M_{\rm {slt}}|X|$,
(2.6), and (3.25), we can verify that the matrices $M_{{\rm {uvec}}}$ and $M_{{\rm
{slvec}}}$ in $Y_{\widetilde{L}}$, and $Y_{\widetilde{U}}$ involved in the bounds given above can be omitted. Thus, the bounds will become concise in form. However, the orders of the matrices in these bounds will increase from ${\nu_1 \times n^2 }$ or ${\nu_2 \times n^2 }$ to ${n^2 \times n^2 }$.}

\end{remark}

\section{PERTURBATION BOUNDS FOR THE QR FACTORIZATION}
Assume that the matrices $A$, $Q$, and $R$ in (1.2) are perturbed as
\begin{eqnarray*}
A\rightarrow A+\Delta A,\ Q\rightarrow Q+\Delta Q,\ R\rightarrow R+\Delta R,
\end{eqnarray*}
where $\Delta A \in {\mathbb{R}^{m \times n}}$, $\Delta Q\in {\mathbb{R}^{m \times n}}$ is such that $(Q+\Delta Q)^T(Q+\Delta Q)=I_n$,
and $\Delta R \in {\mathbb{U}_{n }}$. Thus, the perturbed QR factorization of $A$ is
\begin{eqnarray}
A + \Delta A= (Q + \Delta Q)(R + \Delta R).
\end{eqnarray}
Then
\begin{eqnarray}
(R + \Delta R)^T(R + \Delta R)=(A + \Delta A)^T(A + \Delta A).
\end{eqnarray}
As done in Section 3, here the perturbation $\Delta R$ is also regarded as the unknown matrix.
Expanding (4.2) and considering $A^TA=R^TR$ and (1.2) gives
\begin{eqnarray*}
R^T(\Delta R)+ (\Delta R)^TR=R^T Q^T(\Delta A) +(\Delta A)^TQR+(\Delta
A)^T(\Delta A) -(\Delta R)^T(\Delta R).
\end{eqnarray*}
Left-multiplying the above equation by $R^{-T}$ and right-multiplying it by $R^{-1}$ leads to
\begin{eqnarray*}
(\Delta R)R^{-1}+ R^{-T}(\Delta R)^T=Q^T (\Delta A)R^{-1}
+R^{-T}(\Delta A)^TQ+R^{-T}\left[(\Delta A)^T(\Delta A) -(\Delta
R)^T(\Delta R)\right]R^{-1}.
\end{eqnarray*}
Note that $(\Delta R)R^{-1}$ is upper triangular. Then using the operator `up,' we have
\begin{eqnarray*}
(\Delta R)R^{-1}={\rm {up}}\left[Q^T (\Delta A)R^{-1} +R^{-T}(\Delta
A)^TQ\right]+{\rm {up}}\left(R^{-T}\left[(\Delta A)^T(\Delta A)
-(\Delta R)^T(\Delta R)\right]R^{-1}\right).
\end{eqnarray*}
Applying the operator `vec' to the above equation and using (2.7), (2.4), and (2.8) yields
\begin{eqnarray*}
&&(R^{-T}\otimes I_n){\rm {vec}}(\Delta R)=M_{\rm {up}}\left[R^{-T}\otimes I_n+\left(I_n\otimes R^{-T}\right)\Pi_{nn}\right]{\rm {vec}}\left[Q^T (\Delta A)\right]\\
&&\quad\quad\quad\quad\quad\quad\quad\quad\quad +M_{\rm
{up}}\left(R^{-T}\otimes R^{-T}\right){\rm {vec}}\left[(\Delta
A)^T(\Delta A) -(\Delta R)^T(\Delta R)\right],
\end{eqnarray*}
which together with (2.9) implies
\begin{eqnarray}
&&{\rm {vec}}(\Delta R)=(R^{T}\otimes I_n)M_{\rm {up}}\left[R^{-T}\otimes I_n+\left(I_n\otimes R^{-T}\right)\Pi_{nn}\right]{\rm {vec}}
\left[Q^T (\Delta A)\right]\nonumber\\
&&\quad\quad\quad\quad\quad  +\left(R^{T}\otimes I_n\right)M_{\rm
{up}}\left(R^{-T}\otimes R^{-T}\right){\rm {vec}}\left[(\Delta
A)^T(\Delta A) -(\Delta R)^T(\Delta R)\right].
\end{eqnarray}
Since $\Delta R$ is upper triangular, (2.4), (2.6), and (2.3) together gives
\begin{eqnarray}
{\rm {vec}}(\Delta R)={\rm {vec}}({\rm {ut}}(\Delta R))=M_{{\rm {ut}}}{\rm {vec}}(\Delta R)=M_{{\rm {uvec}}}^TM_{{\rm {uvec}}}{\rm {vec}}(\Delta R)=M_{{\rm {uvec}}}^T{\rm {uvec}}(\Delta R).
\end{eqnarray}
Substituting (4.4) into (4.3) and then premultiplying it by $M_{{\rm {uvec}}}$ and using (2.5), we have
\begin{eqnarray}
&&{\rm {uvec}}(\Delta R)=M_{{\rm {uvec}}}\left(R^{T}\otimes I_n\right)M_{\rm {up}}\left[R^{-T}\otimes I_n+(I_n\otimes R^{-T})\Pi_{nn}\right]
{\rm {vec}}\left[Q^T (\Delta A)\right]\nonumber\\
&&\quad\quad\quad\quad\quad  +M_{{\rm {uvec}}}\left(R^{T}\otimes
I_n\right)M_{\rm {up}}\left(R^{-T}\otimes R^{-T}\right){\rm
{vec}}\left[(\Delta A)^T(\Delta A) -(\Delta R)^T(\Delta R)\right].
\end{eqnarray}
Conversely, left-multiplying (4.5) by $M_{{\rm {uvec}}}^T$ and considering (4.4) and (2.6), we obtain
\begin{eqnarray*}
&&{\rm {vec}}(\Delta R)=M_{{\rm {ut}}}\left(R^{T}\otimes I_n\right)M_{\rm {up}}\left[R^{-T}\otimes I_n+(I_n\otimes R^{-T})\Pi_{nn}\right]
{\rm {vec}}\left(Q^T (\Delta A)\right)\\
&&\quad\quad\quad\quad\quad  +M_{{\rm {ut}}}\left(R^{T}\otimes
I_n\right)M_{\rm {up}}\left(R^{-T}\otimes R^{-T}\right){\rm
{vec}}\left[(\Delta A)^T(\Delta A) -(\Delta R)^T(\Delta R)\right].
\end{eqnarray*}
From the definitions of $M_{{\rm {ut}}}$ and $M_{\rm {up}}$, it is easy to check that
$M_{{\rm {ut}}}\left(R^{T}\otimes I_n\right)M_{\rm
{up}}=\left(R^{T}\otimes I_n\right)M_{\rm {up}}$.
Then the equation (4.3) is equivalent to (4.5).

As a matter of convenience, let
\begin{align*}
G_R=M_{{\rm {uvec}}}\left(R^{T}\otimes I_n\right)M_{\rm
{up}}\left[R^{-T}\otimes I_n+(I_n\otimes R^{-T})\Pi_{nn}\right],
\ H_R=M_{{\rm {uvec}}}\left(R^{T}\otimes I_n\right)M_{\rm
{up}}\left(R^{-T}\otimes R^{-T}\right),
\end{align*}
where $G_R$ is equal to $W_R^{-1}Z_R$ in \cite[Eqn. (3.4.2)]{Chang97a}, but the explicit expression for $W_R^{-1}Z_R$ was not provided in \cite{Chang97a,Chang97}. The fact for equality can be derived from (4.5) and \cite[Eqn. (3.4.2)]{Chang97a} by setting $t=\varepsilon$ in \cite[Eqn. (3.4.2)]{Chang97a} and dropping the higher-order terms. Now, applying the operator `${\rm {uvec^\dag}}$' to (4.5) leads to
\begin{eqnarray}
&&\Delta R={\rm {uvec^\dag}}\Big (G_R{\rm {vec}}\left[Q^T(\Delta
A)\right]+H_R{\rm {vec}}\left[(\Delta A)^T(\Delta A)-(\Delta
R)^T(\Delta R)\right]\Big ).
\end{eqnarray}

In the following, with the help of Lyapunov majorant function and the Banach fixed point principle, we develop the rigorous perturbation bounds for $\Delta R$ based on (4.6).

We first rewrite (4.6) as an operator equation for $\Delta R$,
\begin{eqnarray}
&&\Delta R=\Psi\left(\Delta R,Q^T(\Delta A),\Delta A\right)\\
&&\quad\ \ ={\rm {uvec^\dag}}\Big(G_R{\rm {vec}}\left[Q^T(\Delta
A)\right]+H_R{\rm {vec}}\left[(\Delta A)^T(\Delta A)-(\Delta
R)^T(\Delta R)\right]\Big).\nonumber
\end{eqnarray}
Assuming that $Z\in {\mathbb{U}_{n}}$ and replacing $\Delta R$ in (4.7) with $Z$ leads to
\begin{eqnarray}
Z=\Psi\left(Z,Q^T(\Delta A),\Delta A\right),
\end{eqnarray}
where $\Psi(Z,Q^T(\Delta A),\Delta A)={\rm {uvec^\dag}}\Big(G_R{\rm
{vec}}\left(Q^T(\Delta A)\right)+H_R{\rm
{vec}}\left((\Delta A)^T(\Delta A)-Z^TZ\right)\Big)$.
Let $\left\|Z\right\|_F\leq \rho$ for some $\rho\geq0$,
$\left\|Q^T(\Delta A)\right\|_F= \delta_1$, and $\left\|\Delta
A\right\|_F= \delta_2$. Then, noting (2.1),
\begin{eqnarray*}
\left\|\Psi(Z,Q^T(\Delta A),\Delta
A)\right\|_F\leq\left\|G_R\right\|_2\delta_1+\left\|H_R\right\|_2\delta^2_2+\left\|H_R\right\|_2\rho^2.
\end{eqnarray*}
Thus, setting $\delta=\left[ {\begin{array}{*{20}c}
   \delta_1   \\
   \delta_2  \\
\end{array}} \right]$, we have the Lyapunov majorant function of the operator equation (4.8)
\begin{eqnarray*}
&h(\rho,\delta)=a\delta_1+b\delta^2_2+b\rho^2,
\end{eqnarray*}
where $a=\left\|G_R\right\|_2$ and
$b=\left\|H_R\right\|_2$. Then the Lyapunov majorant equation
is
\begin{eqnarray}
&h(\rho,\delta)=\rho,\ {\textrm{ i.e.}, }\quad
a\delta_1+b\delta^2_2+b\rho^2=\rho.
\end{eqnarray}
Assuming that $\delta\in
\Omega=\{\delta_1\geq0,\delta_2\geq0:
1-4b(a\delta_1+b\delta_2^2)\geq 0\}$,
we have two solutions to the Lyapunov majorant equation (4.9): $\rho_1(\delta)\leq\rho_2(\delta)$ with
\begin{eqnarray*}
\rho_1(\delta):=f_1(\delta):=\frac{2(a\delta_1+b\delta_2^2)}{1+\sqrt{1-4b(a\delta_1+b\delta_2^2)}}.
\end{eqnarray*}
Let the set ${\cal B}(\delta)$ be
\begin{eqnarray*}
{\cal B}(\delta)=\{Z\in {\mathbb{U}_{n}}: \left\|Z\right\|_F\leq
f_1(\delta)\}\subset \mathbb{R}^{n\times n}.
\end{eqnarray*}
It is closed and convex. In this case, the operator $\Psi(\cdot,Q^T(\Delta
A),\Delta A)$ maps the set ${\cal B}(\delta)$ into itself.
Furthermore, when $\delta\in
\Omega_1=\{\delta_{1}\geq0,\delta_2\geq0:
1-4b(a\delta_{1}+b\delta_2^2)> 0\}, $
we have that
the derivative of the function $h(\rho,\delta)$ relative to $\rho$
at $f_1(\delta)$ satisfies
\begin{eqnarray*}
h^{'}_\rho(f_1(\delta),\delta)=1-\sqrt{1-4b(a\delta_1+b\delta_2^2)}\rho<1.
\end{eqnarray*}
Meanwhile, for $Z, \widetilde{Z}\in {\cal B}(\delta)$,
\begin{eqnarray*}
\left\|\Psi(Z,Q^T(\Delta A),\Delta A)-\Psi(\widetilde{Z},Q^T(\Delta
A),\Delta A)\right\|_F\leq
h^{'}_\rho(f_1(\delta),\delta)\left\|Z-\widetilde{Z}\right\|_F.
\end{eqnarray*}
The above facts mean that the operator $\Psi(\cdot,Q^T(\Delta
A),\Delta A)$ is contractive on the set ${\cal B}(\delta)$ when $\delta\in
\Omega_1$. According to the Banach fixed point principle, the
operator equation (4.8) has a
unique solution in the set ${\cal B}(\delta)$ for $\delta\in
\Omega_1$, and so do the operator equation (4.7) and then the matrix equation (4.2). Then
$\left\|\Delta R\right\|_F\leq f_1(\delta)$ for $\delta\in
\Omega_1$. In this case, the unknown matrix $\Delta Q$ in (4.1) is
also determined uniquely.

The above discussions implies another main theorem.
\begin{theorem}
Let the unique QR factorization of $A\in \mathbb{R}^{m\times n}_n$
be as in {\rm (1.2)} and $\Delta A\in \mathbb{R}^{m\times n}$. If
\begin{eqnarray}
\left\|H_R\right\|_2
\left(\left\|G_R\right\|_2\left\|\Delta
A\right\|_F+\left\|H_R\right\|_2\left\|\Delta
A\right\|_F^2\right)<\frac{1}{4},
\end{eqnarray}
then $A+\Delta A$ has the unique QR factorization {\rm (4.1)} and
\begin{eqnarray}
&&\left\|\Delta R\right\|_F\leq \frac{2\left(
\left\|G_R\right\|_2\left\|Q^T(\Delta
A)\right\|_F+\left\|H_R\right\|_2
\left\|\Delta A\right\|_F^2\right)}{1+\sqrt{1-4\left\|H_R\right\|_2 \left(\left\|G_R\right\|_2\left\|Q^T(\Delta A)\right\|_F
+\left\|H_R\right\|_2\left\|\Delta A\right\|_F^2\right)}}\\
&&\quad\quad\quad\ \leq2 \left(\left\|G_R\right\|_2\left\|Q^T(\Delta A)\right\|_F+\left\|H_R\right\|_2\left\|\Delta A\right\|_F^2\right)\\
&&\quad\quad\quad\ <\left(1+2\left\|G_R\right\|_2\right)\left\|\Delta A\right\|_F.
\end{eqnarray}
\end{theorem}
{\raggedleft\em  Proof.}
It is easy to see that the condition (4.10) is more constraining than
the one in $\Omega_1$. Thus, from the  discussions before Theorem
4.1, we derive all results in Theorem 4.1 except the bound (4.13).

After some computations, from (4.10), it follows that
\begin{eqnarray}
2\left\|H_R\right\|_2\left\|\Delta A\right\|_F<{\sqrt{1+\left\|G_R\right\|_2^2}-\left\|G_R\right\|_2}.
\end{eqnarray}
Substituting (4.14) into (4.12) and noting $\left\|Q^T(\Delta
A)\right\|_F\leq\left\|\Delta A\right\|_F$ gives
\begin{eqnarray*}
&&\left\|\Delta R\right\|_F<\left(\sqrt{1+\left\|G_R\right\|_2^2}+\left\|G_R\right\|_2\right)\left\|\Delta A\right\|_F.
\end{eqnarray*}
Using the fact
$\sqrt{1+\left\|G_R\right\|_2^2}\leq1+\left\|G_R\right\|_2$, we
have the bound (4.13).
$\square$

\begin{remark} {\em    According to (4.14), the condition (4.10) can be simplified and strengthened
to
\begin{eqnarray}
\left\|H_R\right\|_2 \left (1+2\left\|G_R\right\|_2
\right)\left\|\Delta A\right\|_F<\frac{1}{2}.
\end{eqnarray}}
\end{remark}
\begin{remark}{\em The following first-order perturbation bound
can be derived from (4.11) or (4.5) by omitting the higher-order
terms
\begin{eqnarray}
\left\|\Delta R\right\|_F\leq
\left\|G_R\right\|_2\left\|Q^T(\Delta A)\right\|_F+{\cal
O}\left(\left\|\Delta A\right\|_F^2\right)
\end{eqnarray}
under the condition
\begin{eqnarray*}
\left\| A^{\dag}\right\|_2\left\|\Delta A\right\|_F<1,
\end{eqnarray*}
which ensures that the unique QR factorization of $A+\Delta A$ exists \cite{Chang97a,
Chang97}.

The bound (4.16) without explicit expression was also derived in \cite{Chang97} by the matrix-vector equation approach,
which is regarded as the optimal first-order bound for
the triangular factor $R$ \cite{Chang97a,Chang97}. }
\end{remark}

\begin{remark} {\em   In \cite{Chang10}, the following rigorous perturbation bound
was derived by the combination of the classic and refined matrix
equation approaches,
\begin{eqnarray}
\left\|\Delta R\right\|_F\leq (\sqrt{6}+\sqrt{3})\left(\mathop {\inf
}\limits_{D \in {\mathbb{D}_{n}}} {{{\sqrt {1 + \zeta _D^2}
}}~{k_2}\left( {{D^{ - 1}}R} \right)} \right)\left\|\Delta
A\right\|_F,
\end{eqnarray}
under the condition
\begin{eqnarray}
\left\| A^{\dag}\right\|_2\left\|\Delta A\right\|_F<\sqrt{3/2}-1.
\end{eqnarray}
In (4.17), $D={\rm diag}(\delta_1,\delta_2,\cdots,\delta_n)$ and
$\zeta _D=\mathop {\max }\limits_{1\leq i < j\leq n} (\delta
_j/\delta _i)$. The discussions in \cite{Chang10} shows that the bound
(4.17) can be much tighter than the previous one derived by the
classic matrix equation approach. From \cite[Eqns. (5.19) and
(5.20)]{Chang97} and the fact $G_R=W_R^{-1}Z_R$ mentioned above, we have
\begin{eqnarray}
1\leq\left\|G_R\right\|_2\leq \mathop {\inf }\limits_{D \in
{\mathbb{D}_{n}}} {{{\sqrt {1 + \zeta _D^2} }}~{k_2}\left({{D^{ -
1}}R} \right)},
\end{eqnarray}
which indicates that the bound (4.13) is tighter than (4.17).

Using the expression of $H_R$ and the definitions of $M_{\rm {uvec}}$ and $M_{\rm {up}}$, we obtain
\begin{eqnarray}
\left\|H_R\right\|_2\geq
\left\| R^{-1}\right\|_2/2=\left\| A^{\dag}\right\|_2/2,
\end{eqnarray}
which together with the first inequality in (4.19) suggests that
\begin{eqnarray*}
\left\|H_R\right\|_2 \left (1+2\left\|G_R\right\|_2
\right)\left\|\Delta A\right\|_F\geq \frac{3}{2}\left\| A^{\dag}\right\|_2\left\|\Delta A\right\|_F.
\end{eqnarray*}
The above inequality is approximately attainable since the
inequality (4.20) and the first inequality in (4.19) are attainable
and approximately attainable \cite{Chang97a,Chang97},
respectively. Moreover, ${1}/{3}>\sqrt{3/2}-1$. So, although the
strengthened condition (4.15) may be more constraining than (4.18),
the former is not so strong. In addition, it should be mentioned
that it is more expensive to estimate the bound (4.13) than that of
(4.17) since the matrix $G_R$ involved in the former contains the
Kronecker products.}
\end{remark}

In the following, we consider the rigorous perturbation bounds for
the triangular factor $R$ of the QR factorization when the
perturbation $\Delta A$ has the form of backward error  resulting
from the standard QR factorization algorithm. That is, $\Delta A\in
\mathbb{R}^{m \times n}$ satisfies (see, e.g., \cite{Higham02,Bai99,Chang01,Zha93}),
\begin{eqnarray}
|\Delta A|\leq \varepsilon C |A|,
\end{eqnarray}
where $C=(c_{ij})\in \mathbb{R}^{m \times m}$, $0\leq c_{ij}\leq 1$, and $\varepsilon\geq 0$ is a small constant. In this case,
\begin{align}
&\left\|\Psi(Z,Q^T(\Delta A),\Delta A)\right\|_F\leq\left\||G_R|{\rm {vec}}\left(|Q^T|C|Q||R|\right)\right\|_F\varepsilon+
\left\||H_R|{\rm {vec}}\left(|R^T||Q^T|C^TC|Q||R|\right)\right\|_F\varepsilon^2\nonumber\\
&\quad\quad\quad\quad\quad\quad\quad\quad\quad\quad \quad +\left\||H_R|\right\|_2\rho^2\nonumber\\
&\quad\quad\quad\quad\quad\quad\quad\quad\quad\ \
\leq\left\||G_R||R^T\otimes
I_n|\right\|_2\left\||Q^T|C|Q|\right\|_F\varepsilon+
\left\||H_R||R^T|\otimes |R^T|\right\|_2\left\| |Q^T|C^TC|Q|\right\|_F\varepsilon^2\nonumber\\
&\quad\quad\quad\quad\quad\quad\quad\quad\quad\quad \quad +\left\||H_R|\right\|_2\rho^2.
\end{align}
From (4.22), we have the Lyapunov majorant function of the
operator equation (4.8) and then (4.7),
\begin{eqnarray*}
h(\rho,\varepsilon)=\widetilde{a}\varepsilon+\widetilde{b}\varepsilon^2+\widetilde{c}\rho^2,
\end{eqnarray*}
where
\begin{eqnarray*}
\widetilde{a}=\left\||G_R||R^T\otimes
I_n|\right\|_2\left\||Q^T|C|Q|\right\|_F,
\end{eqnarray*}
and
\begin{eqnarray*}
\widetilde{b}=\left\||H_R||R^T|\otimes
|R^T|\right\|_2\left\||Q^T|C^TC|Q|\right\|_F,\quad
\widetilde{c}=\left\||H_R|\right\|_2.
\end{eqnarray*}
Then the Lyapunov majorant equation is
\begin{eqnarray*}
h(\rho,\varepsilon)=\rho,\ {\textrm{i.e.}, }\quad
\widetilde{a}\varepsilon+\widetilde{b}\varepsilon^2+{\widetilde{c}}\rho^2=\rho.
\end{eqnarray*}
Similar to the discussions before Theorem 4.1, we have that when $\varepsilon\in \Omega_1$, where
\begin{eqnarray*}
\Omega_1=\left\{\varepsilon\geq0: 1-4\widetilde{c}
(\widetilde{a}\varepsilon+\widetilde{b}\varepsilon^2)>
0\right\},
\end{eqnarray*}
the operator equations (4.8) and (4.7), i.e., the matrix equation (4.2), has a
unique solution in the set
\begin{eqnarray*}
&&{\cal B}(\varepsilon)=\{Z\in {\mathbb{U}_{n}}:
\left\|Z\right\|_F\leq f_1(\varepsilon)\}\subset \mathbb{R}^{n\times
n},
\end{eqnarray*}
where
$f_1(\varepsilon):=\frac{2(\widetilde{a}\varepsilon+\widetilde{b}\varepsilon^2)}
{1+\sqrt{1-4\widetilde{c}(\widetilde{a}\varepsilon+\widetilde{b}\varepsilon^2)}}$.
Then $\left\|\Delta R\right\|_F\leq f_1(\varepsilon)$ for $\varepsilon\in \Omega_1$. In this case, the unknown matrix $\Delta Q$ in (4.1)
is also determined uniquely.

In summary, we have the following theorem.
\begin{theorem}
Let the unique QR factorization of $A\in \mathbb{R}^{m\times n}_n$
be as in {\rm (1.2)} and  $\Delta A\in \mathbb{R}^{m\times n}$ be a
perturbation matrix in $A$ such that {\rm (4.21)} holds. If
\begin{eqnarray}
\widetilde{{c}}(\widetilde{a}\varepsilon+\widetilde{b}\varepsilon^2)<\frac{1}{4},
\end{eqnarray}
then $A+\Delta A$ has the unique QR factorization {\rm (4.1)} and
\begin{align}
&\left\|\Delta R\right\|_F\leq\frac{2(\widetilde{a}\varepsilon+\widetilde{b}\varepsilon^2)}
{1+\sqrt{1-4\widetilde{c}(\widetilde{a}\varepsilon+\widetilde{b}\varepsilon^2)}}\\
&\quad\quad\quad \leq2\left\||G_R||R^T\otimes I_n|\right\|_2\left\||Q^T|C|Q|\right\|_F\varepsilon+
2\left\||H_R||R^T|\otimes |R^T|\right\|_2\left\||Q^T|C^TC|Q|\right\|_F\varepsilon^2 \\
&\quad\quad\quad  < \Big(\left\|
|R|\right\|_2\left\|C|Q|\right\|_F+2\left\||G_R||R^T\otimes
I_n|\right\|_2\left\||Q^T|C|Q|\right\|_F\Big)\varepsilon.
\end{align}
\end{theorem}
{\raggedleft\em  Proof. }  Obviously, we only need to show that the bound (4.26) holds. To see it, we only note the fact
\begin{eqnarray}
0\leq 2\widetilde{b}\varepsilon<{\sqrt{\widetilde{b}/\widetilde{c}+\widetilde{a}^2}-\widetilde{a}}\leq (\widetilde{b}/\widetilde{c})^{1/2}\leq \left\| |R|\right\|_2\left\|C|Q|\right\|_F,
\end{eqnarray}
which can be derived from (4.23) and (2.1).
$\square$

\begin{remark} {\em   Using (4.27), the condition (4.23) can be
simplified and strengthened to
\begin{eqnarray}
\left\||H_R|\right\|_2 \Big(\left\| |
R|\right\|_2\left\|C|Q|\right\|_F+2\left\||G_R||R^T\otimes
I_n|\right\|_2\left\||Q^T|C|Q|\right\|_F\Big)\varepsilon<\frac{1}{2}.
\end{eqnarray}}
\end{remark}

\begin{remark} {\em   From (4.24), we have the following first-order
perturbation bound
\begin{eqnarray}
&&\left\|\Delta R\right\|_F\leq \left\||G_R||R^T\otimes
I_n|\right\|_2\left\||Q^T|C|Q|\right\|_F\varepsilon+{\cal
O}\left(\varepsilon^2\right).
\end{eqnarray}
Replacing $G_R$ with $W_R^{-1}Z_R$ in (4.29) gives the optimal first-order
perturbation bound derived by the matrix-vector equation approach in \cite[Eqn. (8.5)]{Chang01}. In addition, the condition for the bound (4.29) to hold, i.e., for the unique QR factorization $A+\Delta A$ to exist \cite{Chang12a}, is
\begin{eqnarray*}
\left\||R||R^{-1}|\right\|_2\left\|C|Q|\right\|_F\varepsilon<
1.
\end{eqnarray*}
}
\end{remark}

\begin{remark} {\em The following rigorous perturbation bound was
derived by the combination of the classic and refined matrix
equation approaches in \cite{Chang10,Chang12a},
\begin{align}
\left\|\Delta R\right\|_F\leq (\sqrt{6}+\sqrt{3})\left(\mathop {\inf
}\limits_{D \in {\mathbb{D}_{n}}} {{{\sqrt {1 + \zeta _D^2} }}
\left\|{D^{ - 1}}R\right\|_2\left\||R||R^{-1}|D\right\|_2
}\right)\left\|C|Q|\right\|_F\varepsilon,
\end{align}
under the condition
\begin{eqnarray}
\left\||R||R^{-1}|\right\|_2\left\|C|Q|\right\|_F\varepsilon<
\sqrt{3/2}-1.
\end{eqnarray}
It should be claimed that the bound (4.30) is a little different
from the one in \cite{Chang10,Chang12a}. From the discussions in
\cite{Chang12a}, we know that the bound (4.30) can be much smaller than
the one in \cite[Section 6]{Chang01}. Using (2.1), it is seen that
$\left\||Q^T|C|Q|\right\|_F\leq\left\||Q|\right\|_2\left\|C|Q|\right\|_F$.
Meanwhile, from \cite[Eqns. (8.11) and (8.10), and an equation
above (8.7)]{Chang01} and the fact $G_R=W_R^{-1}Z_R$, it follows that
\begin{align}
\left\| |R|\right\|_2\leq\left\|
|G_R||R^T\otimes I_n|\right\|_2\leq
\mathop {\inf }\limits_{D \in {\mathbb{D}_{n}}} {{{\sqrt {1 + \zeta
_D^2} }}\left\|{D^{ -
1}}|R|\right\|_2\left\||R||R^{-1}|D\right\|_2 }.
\end{align}
Thus, when $\left\||Q|\right\|_2= 1$ and
$\left\|{D^{ - 1}}|R|\right\|_2=\left\|{D^{ - 1}}R\right\|_2$, the
bound (4.26) will be tighter than (4.30). As explained in Remark
3.4, a suitable scaling matrix $D$ can make the difference between
$\left\|{D^{ -1}}|R|\right\|_2$ and $\left\|{D^{ - 1}}R\right\|_2$
be unremarkable. See the following examples. So, if
$\left\||Q|\right\|_2= 1$, the bound (4.26) is usually tighter
than (4.30). See Example 4.1 below. Otherwise, since
$\left\||Q^T|C|Q|\right\|_F$ is at most $\left\||Q|\right\|_2$ times
as large as $\left\|C|Q|\right\|_F$, in general, the fact (4.32)
indicates that the bound (4.26) still has advantages. See Example
4.2 below. In addition, we note that the difference between
$\left\||Q^T|C|Q|\right\|_F$ and $\left\|C|Q|\right\|_F$ may
increase as the order $n$ of the involved matrix increases. Example
4.2 given below shows that, in this case, the bound (4.26) still behaves good.

Whereas, the strengthened condition (4.28) may be more constraining than the one (4.31) owing to the first inequality in (4.32) and $\left\||H_R|\right\|_2\geq\left\||R^{-1}|\right\|_2/2$. It is worthy pointing out that the two inequalities mentioned above are attainable \cite{Chang01}. Meanwhile, it is more expensive to estimate the bound (4.26) than that of (4.30), especially when $n$ is large. }
\end{remark}

In the following examples, as done in \cite{Chang01}, we choose the scaling matrix
$D_r$ defined by $D_r={\rm diag}(\left\|R(j,:)\right\|_2)$ and the
scaling matrix $D_e={\rm diag}(\delta_1,\delta_2,\cdots,\delta_n)$
defined as follows: $\delta_1=1/\left\|(D_cR^{-1})(:,1)\right\|_2$;
for $j=2,3,\cdots,n$: $\delta_j=1/\left\|(D_cR^{-1})(:,j)\right\|_2$
if $\left\|(D_cR^{-1})(:,j)\right\|_2\geq
\left\|(D_cR^{-1})(:,j-1)\right\|_2$, otherwise,
$\delta_j=\delta_{j-1}$. Here $D_c={\rm
diag}(\left\|R(j,:)\right\|_1)$. More on methods and explanations of choosing the scaling matrix can be found in \cite{Chang97a} or \cite{Chang98}. In Tables 2--4, we denote
\begin{align*}
&q=\left\||Q^T|C|Q|\right\|_F/\left\|C|Q|\right\|_F,\quad \gamma_{_R}=(\left\| |R|\right\|_2\left\|C|Q|\right\|_F+2\left\||G_R||R^T\otimes I_n|\right\|_2\left\||Q^T|C|Q|\right\|_F)/
\left\|R\right\|_2,\\
&\gamma_{_R}(X)=(\sqrt{6}+\sqrt{3})\left(\sqrt{1 + \zeta _{X}^2}\left\|X^{ - 1}R\right\|_2\left\||R||R^{-1}|X\right\|_2\right)\left\|C|Q|\right\|_F/
\left\|R\right\|_2, \eta_{_{X}}=\left\|{X^{
- 1}}|R|\right\|_2/\left\|{X^{ - 1}}R\right\|_2,
\end{align*}
where $ X=D_r$ or $D_e$, and $t_Y$ the time cost for computing the estimate $Y$. One more statement is that the testing environment is the same as that of Example 3.1.

\begin{example} {\em   This example is from \cite{Chang01}. That is, the test $A$ is the $n\times n$ Kahan matrix:
\begin{align*}
A={\rm diag }(1,s,s^2,\cdots, s^{n-1})
\left[ {\begin{array}{*{20}c}
   1 & { - c} &  \cdots  & { - c}  \\
   {} & 1 &  \cdots  & { - c}  \\
   {} & {} &  \ddots  &  \vdots   \\
   {} & {} & {} & 1  \\
\end{array}} \right],
\end{align*}
where $c=\cos(\theta)$ and $s=\sin(\theta)$. In this case, $R=A$ and $Q=I_n$. Obviously, $\left\||Q|\right\|_2= 1$.
The numerical results for $n=5,10,15,20,25$ with $\theta=\pi/8$ and the corresponding random matrix $C$ produced by the MATLAB function $\mathbf{rand}$ are shown in Table 2, which indicate the expectation claimed in Remark
4.6.

\begin{center}
 \begin{tabular}{|r|r|r|r|r|r|r|r|r|}
    \multicolumn{9}{c}{\centering  Table 2: Comparison of rigorous bounds for the $n\times n$ Kahan matrix } \\\hline
$n$ &$\gamma_{_R}$\quad\quad            &$t_{_{\gamma_{_R}}}$  &$\gamma_{_R}(D_r)$     &$t_{_{\gamma_{_R}(D_r)}}$   &$\eta_{_{D_r}}$ &$\gamma_{_R}(D_e) $  &$t_{_{\gamma_{_R}(D_e)}}$  &$\eta_{_{D_e}}$\\\hline
5  &4.10{\rm e}+01  &0.003 &1.66{\rm e}+02  &0.001 &1.27  &1.79{\rm e}+02  &0.001 &1.05\\
10  &1.48{\rm e}+03 &0.010 &9.00{\rm e}+03  &0.001 &1.27  &1.05{\rm e}+04  &0.001 &1.03\\
15  &4.38{\rm e}+04 &0.036 &3.43{\rm e}+05  &0.002 &1.21  &3.91{\rm e}+05  &0.002 &1.03\\
20  &1.35{\rm e}+06 &0.190 &1.26{\rm e}+07  &0.002 &1.16  &1.40{\rm e}+07  &0.004 &1.03\\
25  &3.87{\rm e}+07 &0.673 &4.15{\rm e}+08  &0.004 &1.13  &4.54{\rm e}+08  &0.004 &1.03\\
\hline
\end{tabular}
\end{center}
}
\end{example}

\begin{example} {\em   Each test matrix has the same form as the one in
Example 3.1. The numerical results for $n=20$, $d_1,d_2\in
\{0.8,1,2\}$, the same random matrix $B$ produced by the MATLAB function $\mathbf{randn}$, and
the same random matrix $C$ produced by the MATLAB function $\mathbf{rand}$ are shown in
Table 3; the numerical results for $n=20, 25, 30, 35, 40, 45, 50, 55$ with $ d_1=d_2=0.8$ and the corresponding random matrices $B$ and $C$ produced by the MATLAB functions $\mathbf{randn}$ and $\mathbf{rand}$, respectively, are shown in Table 4. These results demonstrate the conjectures claimed in Remark
4.6.

\begin{center}
 {\scalebox{1}{\begin{tabular}{|r|r|r|r|r|r|r|r|r|r|r|}
    \multicolumn{11}{c}{\centering  Table 3: Comparison of rigorous bounds for $A=D_1BD_2$ } \\\hline
$d_1$&$d_2$ &$q$ &$\gamma_{_R}$\quad\quad  &$t_{_{\gamma_{_R}}}$&$\gamma_{_R}(D_r)$  &$t_{_{\gamma_{_R}(D_r)}}$ &$\eta_{_{D_r}}$ &$\gamma_{_R}(D_e) $  &$t_{_{\gamma_{_R}(D_e)}}$ &$\eta_{_{D_e}}$ \\\hline
0.8  & 0.8  &2.91   &3.42{\rm e}+02  &0.191 &1.50{\rm e}+03  &0.005 &1.18  &1.45{\rm e}+03  &0.003 &1.00\\
0.8  & 1    &2.91   &9.73{\rm e}+03  &0.192 &5.44{\rm e}+04  &0.003 &1.21  &4.50{\rm e}+04  &0.003 &1.00\\
0.8  & 2    &2.91   &2.29{\rm e}+04  &0.187 &2.90{\rm e}+05  &0.002 &1.07  &1.06{\rm e}+05  &0.003 &1.00\\
1    & 0.8  &3.49   &4.50{\rm e}+02  &0.188 &1.39{\rm e}+03  &0.003 &1.15  &1.32{\rm e}+03  &0.003 &1.00\\
1    &  1   &3.49   &1.52{\rm e}+04  &0.189 &6.62{\rm e}+04  &0.002 &1.32  &4.82{\rm e}+04  &0.003 &1.00\\
1    & 2    &3.49   &2.38{\rm e}+04  &0.190 &6.49{\rm e}+05  &0.003 &1.12  &7.56{\rm e}+04  &0.003 &1.00\\
2    & 0.8  &2.00   &4.38{\rm e}+02  &0.187 &3.77{\rm e}+03  &0.003 &1.15  &3.11{\rm e}+03  &0.003 &1.02\\
2    & 1    &2.00   &3.39{\rm e}+02  &0.191 &1.37{\rm e}+05  &0.003 &1.17  &2.33{\rm e}+04  &0.006 &1.03\\
2    & 2    &2.00   &8.02{\rm e}+03  &0.188 &1.94{\rm e}+06  &0.003 &1.05  &5.48{\rm e}+04  &0.002 &1.00\\
\hline
\end{tabular}}}
\end{center}

\begin{center}
  {\scalebox{1}{\begin{tabular}{|r|r|r|r|r|r|r|r|r|r|}
    \multicolumn{10}{c}{\centering  Table 4: Comparison of rigorous bounds for $A=D_1BD_2$ with $d_1=d_2=0.8$ } \\\hline
$n$&$q$ &$\gamma_{_R}$\quad\quad &$t_{_{\gamma_{_R}}}$ &$\gamma_{_R}(D_r)$  &$t_{_{\gamma_{_R}(D_r)}}$ &$\eta_{_{D_r}}$ &$\gamma_{_R}(D_e) $ &$t_{_{\gamma_{_R}(D_e)}}$ &$\eta_{_{D_e}}$  \\\hline
20    &2.99   &4.56{\rm e}+02  &0.190 &1.24{\rm e}+03  &0.007 &1.19  &1.51{\rm e}+03  &0.009 &1.08\\
25    &3.22   &8.42{\rm e}+02  &0.662 &1.96{\rm e}+03  &0.005 &1.10  &2.39{\rm e}+03  &0.005 &1.20\\
30    &3.33   &7.64{\rm e}+02  &1.914 &2.93{\rm e}+03  &0.007 &1.27  &3.26{\rm e}+03  &0.006 &1.05\\
35    &3.31   &7.29{\rm e}+02  &4.688 &1.68{\rm e}+03  &0.008 &1.22  &3.05{\rm e}+03  &0.008 &1.06\\
40    &3.34   &1.11{\rm e}+03  &10.69 &3.06{\rm e}+03  &0.011 &1.15  &4.50{\rm e}+03  &0.009 &1.14\\
45    &3.45   &1.04{\rm e}+03  &21.35 &3.48{\rm e}+03  &0.013 &1.18  &4.69{\rm e}+03  &0.012 &1.07\\
50    &3.50   &7.33{\rm e}+02  &39.81 &2.65{\rm e}+03  &0.012 &1.12  &4.31{\rm e}+03  &0.012 &1.00\\
55    &3.46   &1.51{\rm e}+03  &69.81 &3.93{\rm e}+03  &0.014 &1.28  &6.75{\rm e}+03  &0.014 &1.13\\
\hline
\end{tabular}}}
\end{center}
}
\end{example}

\begin{remark}{\em
As done in the proof of Theorem 3.1 and Remark 3.5 and noting the fact $M_{{\rm {ut}}}(R^{T}\otimes I_n)M_{\rm {up}}=(R^{T}\otimes I_n)M_{\rm {up}}$, we can check that the matrix $M_{{\rm {uvec}}}$ in $G_R$ and $H_R$ involved in the bounds given in this section can be omitted. In this case, the forms of these bounds will become concise, however, the orders of the matrices in these bounds will increase.}

\end{remark}

\section{CONCLUDING REMARKS}
In this paper, we propose a new approach to present the rigorous
perturbation analysis for the LU and QR factorizations, and obtain
new rigorous perturbation bounds with explicit expressions, which improve the previous
ones in \cite{Chang10} and \cite{Chang12a}. As the special case, the optimal first-order
perturbation bounds with explicit expressions for the two factorizations are also given. The new approach can also be used to derive the
rigorous perturbation bounds for the Cholesky factorization and the
Cholesky downdating problem \cite{Chang96,Chang97a,Chang98b}. The derived bounds for the
Cholesky factorization are the same as the ones in \cite{Chang96,Chang97a} obtained by
the combination of the matrix-vector equation approach and Theorem
3.1 in \cite{Stewart73}, but have the explicit expressions. Actually, noting the conditions and proof of
Theorem 3.1 in \cite{Stewart73}, we find that the approach in \cite{Chang96,Chang97a} can be
regarded as a special case of the approach in this paper. Furthermore, the new approach can also be generalized to apply the block matrix factorizations such as the block LU, SR, and Cholesky-like factorizations \cite{Be}.

Although the explicit expressions of the new rigorous perturbation
bounds and the optimal first-order perturbation bounds are provided,
it is still expensive to estimate these bounds directly as the
spectral norm of the large sparse matrices is involved. To reduce
the computational cost, we can use the fact that, for any matrix
$X$, $\left\| X\right\|_2^2\leq\left\| X\right\|_1\left\|
X\right\|_\infty$. However, in this case, the bounds will be
weakened. In addition, some techniques on sparse matrix (see
e.g., \cite{Davis06}) may be used to overcome the above difficulties. We
will consider this topic in the near future.




\end{document}